\RequirePackage{ifpdf}

\RequirePackage{fix-cm}

\documentclass[%
floatfix,
showkeys,
twocolumn, %
nofootinbib, %
superscriptaddress, %
]{revtex4-1}

\ifpdf
  \usepackage{cmap}
\fi
\usepackage{textcomp}

\usepackage[T1,T2A]{fontenc}

\usepackage[utf8]{inputenx}
\input{ix-utf8enc.dfu}

\usepackage[german,russian,english]{babel}

\usepackage{amsmath}
\usepackage{amssymb}

\usepackage{mathtools}
\mathtoolsset{
showonlyrefs,
mathic = true
}

\allowdisplaybreaks

\usepackage{hyperref}
\hypersetup{backref,
 colorlinks=false}
\hypersetup{pdfborder=0 0 0}

\usepackage{microtype}
\UseMicrotypeSet[protrusion]{alltext}

\usepackage{graphicx}

\usepackage[scanall]{psfrag}

\usepackage{listings}
\usepackage{listingsutf8}
\usepackage[frozencache]{minted}

\setminted{breaklines=true}

\usepackage{fixltx2e}

\usepackage{nicefrac}

\makeatletter
\def\ps@pprintTitle{%
     \let\@oddhead\@empty
     \let\@evenhead\@empty
     \let\@oddfoot\@empty
     \let\@evenfoot\@oddfoot}
\makeatother

\usepackage{physics}

\begin{document}

\graphicspath{{image/symbolic-rk/en/}{image/symbolic-rk/}{image/}}

\title{Using a template engine as a computer algebra tool}

\author{Migran N. Gevorkyan}
\email{gevorkyan-mn@rudn.ru}
\affiliation{Department of Applied Probability and Informatics\\
  Peoples' Friendship University of Russia (RUDN University)\\
  6, Miklukho-Maklaya St., Moscow, 117198, Russian Federation}

\author{Anna V. Korolkova}
\email{korolkova-av@rudn.ru}
\affiliation{Department of Applied Probability and Informatics\\
  Peoples' Friendship University of Russia (RUDN University)\\
  6, Miklukho-Maklaya St., Moscow, 117198, Russian Federation}

\author{Dmitry S. Kulyabov}
\email{kulyabov-ds@rudn.ru}
\affiliation{Department of Applied Probability and Informatics\\
  Peoples' Friendship University of Russia (RUDN University)\\
  6, Miklukho-Maklaya St., Moscow, 117198, Russian Federation}
\affiliation{Laboratory of Information Technologies\\
  Joint Institute for Nuclear Research\\
  6, Joliot-Curie, Dubna, Moscow region, 141980, Russian Federation}

\begin{abstract}
In research problems that involve the use of numerical methods for
solving systems of ordinary differential equations (ODEs), it is
often required to select the most efficient method for a particular
problem. To solve a Cauchy problem for a system of ODEs, Runge--Kutta
methods (explicit or implicit ones, with or without step-size
control, etc.) are employed. In that case, it is required to search
through many implementations of the numerical method and select
coefficients or other parameters of its numerical scheme. This paper
proposes a library and scripts for automated generation of routine
functions in the Julia programming language for a set of numerical
schemes of Runge–Kutta methods. For symbolic manipulations, we use a
template substitution tool. The proposed approach to automated
generation of program code allows us to use a single template for
editing, instead of modifying each individual function to be
compared. On the one hand, this provides universality in the
implementation of a numerical scheme and, on the other hand, makes
it possible to minimize the number of errors in the process of
modifying the compared implementations of the numerical method. We
consider Runge--Kutta methods without step-size control, embedded
methods with step-size control, and Rosenbrock methods with
step-size control. The program codes for the numerical schemes,
which are generated automatically using the proposed library, are
tested by numerical solution of several well-known problems.
\end{abstract}


\maketitle

\section{Introduction}
\label{sec:intro}

  Runge–Kutta methods are basic numerical methods for solving nonrigid
  systems of ordinary differential equations (ODEs). Numerical schemes
  of high orders (tenth and higher) with step-size control and dense
  output are well known. For many programming languages, there are
  libraries that implement the most efficient numerical schemes. The
  most well debugged codes were generated for Fortran 77 by Ernst
  Hairer and his colleagues; the codes were described in~\cite{L_Hairer_1,L_Hairer_2} and
  are available online~\cite{L_FortranCodes}. These routines implement the methods
  DOPRI5~\cite{L_DORPRI1980} and DOPRI853~\cite{L_DORPRI1981}, which feature step-size control and
  dense output. The latter method also supports order switching
  between $8$ and $5$ for better performance when solving smooth problems.

  The DOPRI5 and DOPRI853 subroutines are very well optimized and are
  included in many libraries and math packages, e.g., Matlab~\cite{L_Matlab},
  Octave~\cite{L_Octave}, SciPy~\cite{L_Scipy}, SciLab~\cite{L_Scilab}, Boost C++~\cite{L_Boost}, etc. However, for
  educational and research purposes, it is often required to search
  through a number of Runge–Kutta methods to select the optimal one
  for a particular problem~\cite{kulyabov:2016:rk-stochastic,kulyabov:2017:ceur:vol2064:runge-kutta,kulyabov:2018:dccn:rk-stochastic}. In addition, the efficiency of a
  particular numerical scheme also significantly depends on a problem
  to be solved.

  To reliably compare different numerical schemes, we need a universal
  code that implements embedded methods for any set of
  coefficients. In this case, all implementations must be unified and
  differ only in their sets of coefficients.

  To solve this problem, we propose a library for Julia~\cite{L_JuliaLang_1,L_JuliaLang_2,L_JuliaLang_3} that
  consists of two parts: a computation part that implements algorithms
  of the Runge–Kutta type and a construction part that uses symbolic
  computations to generate specific versions of the Runge–Kutta
  algorithm and features automatic generation of function code from
  several ready-made templates. Initially, we considered implementing
  the construction part by using a universal computer algebra system,
  e.g., SymPy~\cite{kulyabov:2020:rng-cas-sympy::en}. However, in the process of implementation, it
  turned out that the main requirement for the computer algebra system
  was the support of template comparison. That is why we decided to
  use a lighter tool for symbolic manipulations. More specifically,
  instead of a full-featured computer algebra system, we use a
  template substitution tool. This approach, although somewhat exotic,
  is justified by a higher performance of the functions generated. For
  code generation, we use Python~\cite{L_Python} and Jinja2 template engine~\cite{L_Jinja2}.

\subsection{Paper  structure}
\label{sec:structure}

  Section~\ref{sec:rk} briefly discusses embedded Runge--Kutta methods and
  step-size control strategy. Section~\ref{sec:program} describes the developed
  library for Julia. Section~\ref{sec:test} presents the results of testing the
  functions generated automatically using the proposed library for the
  numerical solution of systems of ODEs. In particular, that section
  considers in detail the solution of the restricted three-body
  problem (problem of computing Arenstorf orbits).

\section{Runge--Kutta methods}
\label{sec:rk}

In this paper, the Runge–Kutta methods are applied to a Cauchy problem
formulated for a system of ODEs. The theory of Runge--Kutta methods is
well known and described in detail in~\cite{L_Hairer_1, L_Hairer_2, L_Butcher}. Thus, in this
section, we confine ourselves only to the basic formulas while briefly
discussing the step-size control strategy implemented in our library. 

\subsection{Problem statement}

  Suppose that we have two smooth functions:
  $y^{\alpha}(t)\colon [t_0, T] \to \mathbb{R}^{N}$ an unknown
  function and a known function
  $f^{\alpha}(t, y^{\beta}(t))\colon \mathbb{R}\times\mathbb{R}^N \to
  \mathbb{R}^{N}$, where $[t_0, T] \in \mathbb{R}$ and
  $\alpha=1,\ldots,N$. Suppose also that the value of the function at
  the initial instant is known: $y^{\alpha}_0 =
  y^{\alpha}(t_0)$. Then, the Cauchy problem for a system of ODEs is
  formulated as follows:
\begin{equation}
  \label{Cauchy_ODE_0}
  \left\{
  \begin{aligned}
    &\dot{y}^{\alpha}(t) = f^{\alpha}(t, y^{\beta}(t)),\\
    & y^{\alpha}(t_0) = y^{\alpha}_0,\;\alpha,\beta = 1,\ldots,N.\\
  \end{aligned}
  \right.
\end{equation}

System~\eqref{Cauchy_ODE_0} can be written component by component: 
\begin{equation}
  \label{Cauchy_ODE}
 \left\{
 \begin{aligned}
   &\dot{y}^1(t) = f^1(t, y^1(t),\ldots,y^N(t)),\\
  %  &\dot{y}^2(t) = f^2(t, y^1(t),\ldots,y^N(t)),\\
   &\vdots\\
   &\dot{y}^N(t) = f^N(t, y^1(t),\ldots,y^N(t)),\\
   & y^{\alpha}(t_0) = y^{\alpha}_0.\\ 
 \end{aligned}
 \right.
\end{equation}

  On a segment $[t_0, T]$, we define a grid from a set of points
  $t_0 < t_1 < t_2 < \ldots < t_k < \ldots < t_N = T$ with a grid step
  $h_{k+1} = t_{k+1} - t_{k}$. Based on a certain rule called a
  \emph{numerical scheme}, each point of the grid is associated with a
  certain value $y^{\alpha}_k$ that approximates (with desired
  accuracy) the solution of the system of ODEs at this point, i.e.,
  $y^{\alpha}_k \approx y^{\alpha}(t_k)$. To estimate approximation
  error, we use the norm $\| y^{\alpha}(t_{k}) - y^{\alpha}_{k} \|$.

\subsection{Embedded explicit Runge--Kutta methods}

The error of the Runge--Kutta method can also be estimated in other
way. The idea is that, in addition to the main solution $y^{\alpha}_{m}$ at a certain
point, we also consider an auxiliary solution $\Hat{y}^{\alpha}_{m}$ found by a Runge--Kutta
method of an adjacent order. The difference between these solutions
serves as an estimate for the local error of the lower-order
method. This estimate facilitates the selection of a variable
integration step. The methods that use the local error estimate are
called \emph{embedded} Runge--Kutta methods~\cite{L_Hairer_1}. In this case, the method $y^{\alpha}_{m}$ is
called \emph{main}, while the method $\Hat{y}^{\alpha}_{m}$ is called \emph{embedded}.

An embedded explicit Runge–Kutta method for Cauchy problem~\eqref{Cauchy_ODE} is given by the following formulas: 
\begin{equation*}
  \begin{aligned}
    & k^{1\alpha} = f^{\alpha}(t_{m}, y^{\beta}_{m}),
    \\
    &
    \begin{multlined}
      k^{i\alpha} = f^{\alpha}\big(t_{m} + c^{i}h_{m}, y^{\alpha}_{m}
      + h\sum^{i-1}_{j=1}a^{i}_{j}k^{j \beta}\big),
      \\
      i=2,\ldots,s,
    \end{multlined}
    \\
    &
    \begin{multlined}
      y^{\alpha}_{m+1}= y^{\alpha}_{m} + h_{m}(b_{1}k^{1\alpha} +
      b_{2}k^{2\alpha}
      + {} \\ {} +
      \ldots + b_{s-1}k^{s-1,\alpha} +
      b_{s}k^{s\alpha}),
    \end{multlined}
    \\
    &
    \begin{multlined}
      \hat{y}^{\alpha}_{m+1}= y^{\alpha}_{m} +
      h_{m}(\hat{b}_{1}k^{1\alpha} + \hat{b}_{2}k^{2\alpha}
      + {} \\ {} +
      \ldots +
      \hat{b}_{s-1}k^{s-1,\alpha} + \hat{b}_{s}k^{s\alpha}),
    \end{multlined}
  \end{aligned}
\end{equation*}
where $N$ is the number of equations in the system of ODEs and $s$ is the
number of stages in the numerical method. The Latin indices $i,j,l=1,\ldots,s$ are
associated with the numerical scheme, the Greek indices $\alpha$ and $\beta$ are
associated with the system of ODEs, while the subscript $m$ indicates the
number of a step.

The orders of approximation by the grid functions $y^{\alpha}_{m}$ and $\hat{y}^{\alpha}_{m}$ are different,
which makes it possible to compute approximation error at each
step. Thus, to specify the order of the embedded method, we use the
notation $p(\Hat{p})$, where $p$ and $\Hat{p}$ are the orders of the main and embedded
Runge--Kutta methods, respectively.

The values $c^{i},a^{i}_{j},b_{j},\hat{b}_{j}$ completely define the numerical scheme and are called
\emph{coefficients of the method}; they are usually grouped in a Butcher
table (named after John C. Butcher):
\begin{equation*}
  \begin{array}{c|cccccc}
  0  & 0 & 0 &\ldots & 0 & 0\\
  c^2  & a^{2}_1 & 0 &\ldots & 0 & 0\\
  c^3  & a^{3}_1 &a^{3}_2 &\ldots &0 & 0\\
  \vdots  &\vdots  &\vdots &\ddots &\vdots &\vdots\\
  c^{s-1} & a^{s-1}_1 &a^{s-1}_2 &\ldots & 0 & 0\\
  c^{s}  & a^{s}_1 &a^{s}_2 &\ldots &a^{s}_{s-1} & 0\\ \hline
       &b_{1} &b_{2} &\ldots &b_{s-1} & b_{s}\\
       &\hat{b}_{1} &\hat{b}_{2} &\ldots &\hat{b}_{s-1} & \hat{b}_{s}
  \end{array}
\end{equation*}

Finding the coefficients of the Runge--Kutta method for orders higher
than fourth is a separate complex problem, which falls outside the
scope of this paper and is described in detail in~\cite{L_Hairer_1}. It should be
noted, however, that it is common practice to impose the following
additional conditions on the coefficients $c^{i}$ and $a^{i}_{j}$: $c^{i} = a^{i}_{1} + \ldots + a^{i}_{s}$.

  To date, the coefficients for many embedded Runge–Kutta methods up
  to the tenth order of accuracy have been found. The most popular
  methods are those found by J.R. Dorman and
  P.J. Prince~\cite{L_DORPRI1980,L_DORPRI1981,L_DORPRI1986},
  E. Fehlberg~\cite{fehlberg:1969:klassische-runge-kutta-funfter,fehlberg:1970:klassische-runge-kutta-vierter},
  and J.R. Cash and A.H. Karp~\cite{L_CashKarp1990}.

\subsection{Step-size control strategy}

There are various strategies for step-size (h) control depending on
the local error. The choice of a particular strategy generally depends
on a problem to be solved. Here, we employ the classical algorithm
proposed by E. Hairer\cite{L_Hairer_2}, which is based on ideas borrowed from
control theory and works well for most nonrigid problems.

Suppose that $y^{\alpha}_{m}$ and $\hat{y}^{\alpha}_{m}$ are two numerical solutions of different
approximation orders that are computed at the step $m$ of the
algorithm. It is required that , where the desired scale $sc^{\alpha}$ is found
by the formula
\begin{equation*}
  sc^{\alpha} = A_{tol} + \max(|y^{\alpha}_{m}|, |\hat{y}^{\alpha}_{m}|)R_{tol},
\end{equation*}
where $A_{tol}$ and $R_{tol}$ are the desired absolute and relative tolerances,
respectively. The permissible error at the step $m$ is found as a root
mean square:
\begin{equation*}
  E_{m} = \sqrt{\dfrac{1}{N} \sum\limits_{\alpha = 1}^{N}{\left(\dfrac{y^{\alpha}_{m} - y^{\alpha}_{m}}{sc^{\alpha}}\right)^{2}}}.
\end{equation*}

A new step size is computed by the formula 
\begin{equation*}
  h_{m+1} = h_{m}/\max\left(f_{\min}, \min\left(f_{\max}, E^{a}_{m}E^{-b}_{m-1}/f_{s}\right)\right),
\end{equation*}
  where $E_{m-1}$ is the error found at the previous step. The factor
  $f_{s}$ (used to prevent a sharp increase in the step size) is
  generally $0.9$ or $0.8$. The exponents $a$ and $b$ are selected
  depending on a problem. In~\cite{L_Hairer_2}, the following
  universal values were recommended: $a = 0.7/p - 0.75b$ and
  $b = 0.4/p$, where $p$ is the order of approximation. The factors
  $f_{\min}$ and $f_{\max}$ set the step-size variation boundaries and
  also depend on a problem to be solved. In practice, the values
  $f_{\min} = 0.1$ and $f_{\max} = 5.0$ are generally selected.

  In the case of $E_{m} < 1$, the computed values $y^{\alpha}_{m}$ are
  considered satisfactory and the method goes to the next step. If
  $E_{m} > 1$, then the results are considered unsatisfactory and the
  current step is repeated with a new value of $h_{m}$:
\begin{equation*}
  h_{m} \leftarrow h_{m} / \min\left(f_{\max}, E^{a}_{m}/f_{s}\right).
\end{equation*}

At the next step $m+1$, as the initial value $y^{\alpha}_{m}$, we can use either
directly $y^{\alpha}_{m}$ or the value yielded by the embedded method $\hat{y}^{\alpha}_{m}$.

\section{Description of the routine}
\label{sec:program}

\subsection{Motivation to use code generation}

Universal implementation of an embedded explicit Runge–Kutta method
implies that, as parameters, the routine must receive arrays of
coefficients from Butcher tables, which are then used in
computations. Arrays seem to be a natural way to store these
coefficients. However, $a^{i}_{j}$ is a lower-diagonal matrix and, when storing
it as an $s\times s$ array, more than half of the memory allocated for the
array is spent on storing zeros. The coefficients $c^{i}$,
$b_{j}$, and $\hat{b}_{j}$ also
often contain zeros, storing which is unreasonable.

That is why most codes that implement embedded explicit Runge–Kutta
methods use a set of named constants (rather than arrays) to store the
coefficients. Moreover, in most languages, scalar operations are
faster than array operations.

In Julia, standard arrays are dynamic, which is why the overhead of
storing a two-dimensional $s \times s$ array exceeds the overhead of storing
$s \times s$ named constants.

To preserve the requirement for universality of the generated code
while speeding up computations and reducing memory consumption, we
decided to use automatic code generation based on one template for
each individual method.

In addition to gaining in performance, automatic code generation
allows us to add or modify all functions at once (rather than each
function individually) by editing only one template. This also makes
it possible to reduce the number of errors and generate different
variants of functions for different purposes.

\subsection{Description of the Runge--Kutta Generator}

As a language for code generation, we chose Python because it supports
various means for text manipulation. In addition, the standard Python
library contains the \verb|Fraction| data type, which allows us to set the
Runge--Kutta coefficients as rational fractions and then convert them
to real-valued form with desired order of accuracy.

In addition to Python, we also employ the Jinja2 template processing
library~\cite{L_Jinja2}. This template engine was originally designed for
generating HTML pages. However, it has a very flexible syntax and can
be used as a universal tool for generating text files of any kind,
including source codes in any programming language. In combination
with Jinja2, we use the \verb|numpy| library~\cite{L_Scipy} to process arrays of
coefficients.

The templates for generating functions are stored in the files
\verb|rk_method.jl| and \verb|erk_method.jl|, which are Julia source codes with
Jinja2 instructions. The template engine allows us to place the entire
code generation logic in a template and input only the data associated
with the method.

Information about numerical schemes is stored in a JSON file, where
each method is a JSON object similar to the following one:
\begin{minted}{json}
{"name": "Name (one word)",
"description": "Brief description",
"stage": ,
"order": ,
"extrapolation_order": ,
"a": [["0" , "0" , "0" ], ["1" , "0" ,
"0" ], ["1/4" , "1/4" , "0" ]],
"b": ["1/2" , "1/2" , "0" ],
"b_hat": ["1/6" , "1/6" , "2/3" ],
"c": ["0" , "1" , "1/2" ]}
\end{minted}

The name of the JSON object is used as a name of a subsequently
generated function. The arrays \verb|a|, \verb|b|, \verb|b_hat|, and \verb|c| can have either
numeric or string type. If the method's coefficients are set as
rational fractions, then we can define them as \verb|"m/n"|. Further, they
are converted to a \verb|Fraction| object of the Python standard library; in
the body of a function generated, they are represented as
double-precision decimal fractions with $17$ significant digits.

The list of these objects is sequentially processed by the scripts
\verb|erk_generator.py| and \verb|rk_generator.py|. By using the information about
the methods, the scripts generate codes of several functions (in
Julia) for each method.

These scripts generate functions for 16 embedded Runge--Kutta
methods. Low-order methods are borrowed from~\cite{L_Hairer_1}. For the methods of
the fifth order and higher, we use the coefficients from~\cite{L_DORPRI1980,L_DORPRI1981,L_DORPRI1986}
and~\cite{fehlberg:1969:klassische-runge-kutta-funfter,fehlberg:1970:klassische-runge-kutta-vierter}. The method from~\cite{L_CashKarp1990} is implemented separately as it
supports order switching and has a complex step-size control
algorithm. To add user methods, it is required to create a JSON object
(see above) and run the scripts to generate the corresponding
functions for these methods.

\subsection{Description of the functions generated}

The developed library has a structure typical for Julia modules. The
entire source code is stored in the \verb|src| directory. The subdirectory
\verb|src/generated| contains files with automatically generated functions,
which, in turn, are included in the main module file \verb|src/RungeKutta.jl|
by using the \verb|include| directive. The code responsible for step-size
control is stored in the file \verb|StepControl.jl| and is common to all
methods.

For each of the embedded Runge--Kutta methods, three functions are
generated, two of which have the following form:
\begin{minted}{julia}
  ERK(func, A_tol, R_tol, x_0, t_start, t_stop) -> (T, X)
  ERK(func, A_tol, R_tol, x_0, t_start, t_stop, last) -> (t, x)
\end{minted}

The functions have the same name and their lists of arguments differ
only in the last argument. With Julia supporting multiple dispatch
(function overloading), the compiler by itself determines which
implementation needs to be called in each case.

In the functions described above,
\begin{itemize}  
\item \verb|func::Function| is the right-hand side of the system of
  ODEs $\dot{x}^{\alpha}(t) = f^{\alpha}(t, x^{\beta})$:
  \verb|func(t,x)|,
  where \verb|t::Float64| is time and \verb|x::VectorFloat64| is
  the value of the function $x(t)$; the argument $t$ must always be
  included in the function call, even if the system is autonomous and
  obviously does not depend on time;
  
\item the arguments \verb|A_tol| and \verb|R_tol| are the absolute and
  relative accuracies of the methods, respectively; they must have the
  type \verb|Float64|;
  
\item \verb|x_0::Vector{Float64}| is the initial value of the function
  $x^{\alpha}_{0} = x^{\alpha}(t_0)$;
  
\item \verb|t_start| and \verb|t_stop| are the start and end points of
  the integration interval, respectively; they have the type
  \verb|Float64|;
  
\item if the argument \verb|last::Bool| is specified, then the end
  point of the integration interval $t_{n}$ and vector $x^{\alpha}_{n}
  \approx x^{\alpha}(t_n)$ are returned; otherwise, the arrays
  \verb|T::Vector{Float64}| and \verb|Matrix{Float64}| are returned.  
\end{itemize}

For the Runge--Kutta methods without step-size control, similar
functions are generated with their only difference being the absence
of the arguments \verb|A_tol| and \verb|R_tol|. Instead of them, the argument \verb|h|,
which specifies a step size used for computation, is passed.

To implement the step-size control algorithm, another function is generated:
\begin{minted}{julia}
  ERK_info(func, A_tol, R_tol, x_0, t_start, t_stop) -> (accepted_t, accepted_h, rejected_t, rejected_h)
\end{minted}
This function returns the following arrays:
\begin{itemize}  
\item \verb|accepted_t|: all grid points at which the computation error is considered satisfactory;
\item \verb|accepted_h|: all accepted step sizes;
\item \verb|rejected_t|: all grid points at which the computation error is considered unsatisfactory;
\item \verb|rejected_h|: all rejected step sizes;
\item \verb|errors|: values of the local error $E_{n}$. 
\end{itemize}

All returned values have the type \verb|Vector{Float64}|.
  \section{Testing the codes generated for the numerical schemes}
\label{sec:test}

To test the methods, we use three systems of differential equations
considered in~\cite{L_Hairer_2}. The first system is a system of van
der Pol equations
\begin{gather*}
  \dv[2]{x}{t} - \mu(1-x^2) \dv{x}{t} + x = 0,
  \\
  \left\{
  \begin{aligned}
      & \Dot{x}_{1}(t) = x_{2}(t),\\
      & \Dot{x}_{2}(t) = \mu (1-x^{2}_{1}(t)) x_{2}(t) - x_{1}(t),
  \end{aligned}
  \right.
\end{gather*}
with the initial values
\begin{equation*}
  \mu = 1,
  \quad
  \mathbf{x} = (0, \sqrt{3})^T,
  \quad
  0 \leqslant t \leqslant 12,
\end{equation*}
where the coefficient $\mu$ characterizes the nonlinearity and decay of oscillations.

The second system is a system of equations for a rigid body without
external forces, i.e., a system of Euler's equations for a rigid body:
\begin{equation*}
  \left\{
  \begin{aligned}
      & \dot{x}_{1}(t) = I_{1} x_{2}(t) x_{3}(t),\\
      & \dot{x}_{2}(t) = I_{2} x_{1}(t) x_{3}(t),\\
      & \dot{x}_{3}(t) = I_{3} x_{1}(t) x_{2}(t),\\
      & I_1 = -2,\;I_2 = 1.25,\; I_3 = -0.5,\\
      & \mathbf{x} = (0, 1, 1)^T,\; t = [0, 12].
  \end{aligned}
  \right.
\end{equation*}

To test the step-size control algorithm, we apply the numerical scheme
of the embedded Runge--Kutta method with $p=4$ and $\Hat{p}=3$ to a system of
Brusselator equations
\begin{equation}
  \label{eq:01}
  \left\{
  \begin{aligned}
      & \dot{x}_{1}(t) = 1 + x^{2}_{1}x_{2} - 4x_{1},\\
      & \dot{x}_{2}(t) = 3x_{1} - x^{2}_{1} x_{2}.
  \end{aligned}
  \right.
\end{equation}

\begin{figure}[!h]
  \centering
  \includegraphics[width=1.0\linewidth]{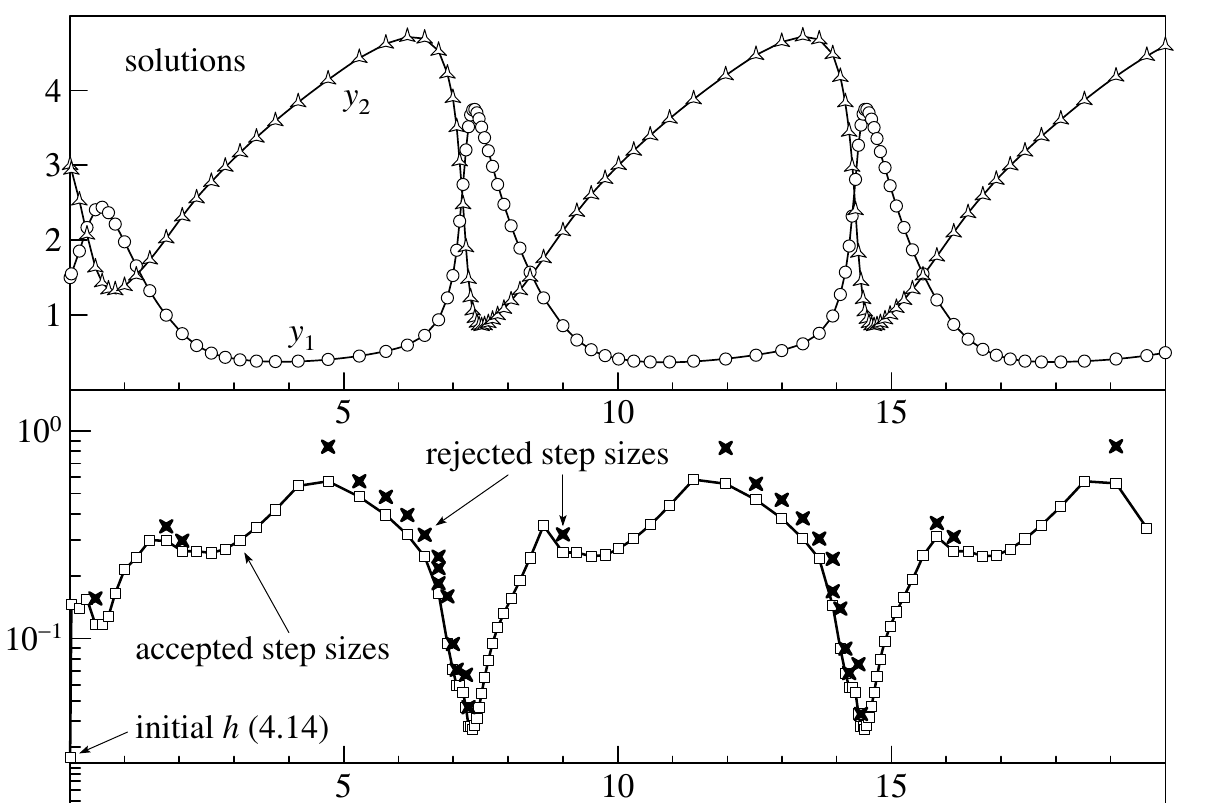}
\caption{Solution of system~\eqref{eq:01} with rejected and accepted steps ~\cite[p.~170, Fig.~4.1]{L_Hairer_1}}
  \label{fig:erk43b:hairer}
\end{figure}

\begin{figure}[!h]
  \centering
  \includegraphics[width=1.0\linewidth]{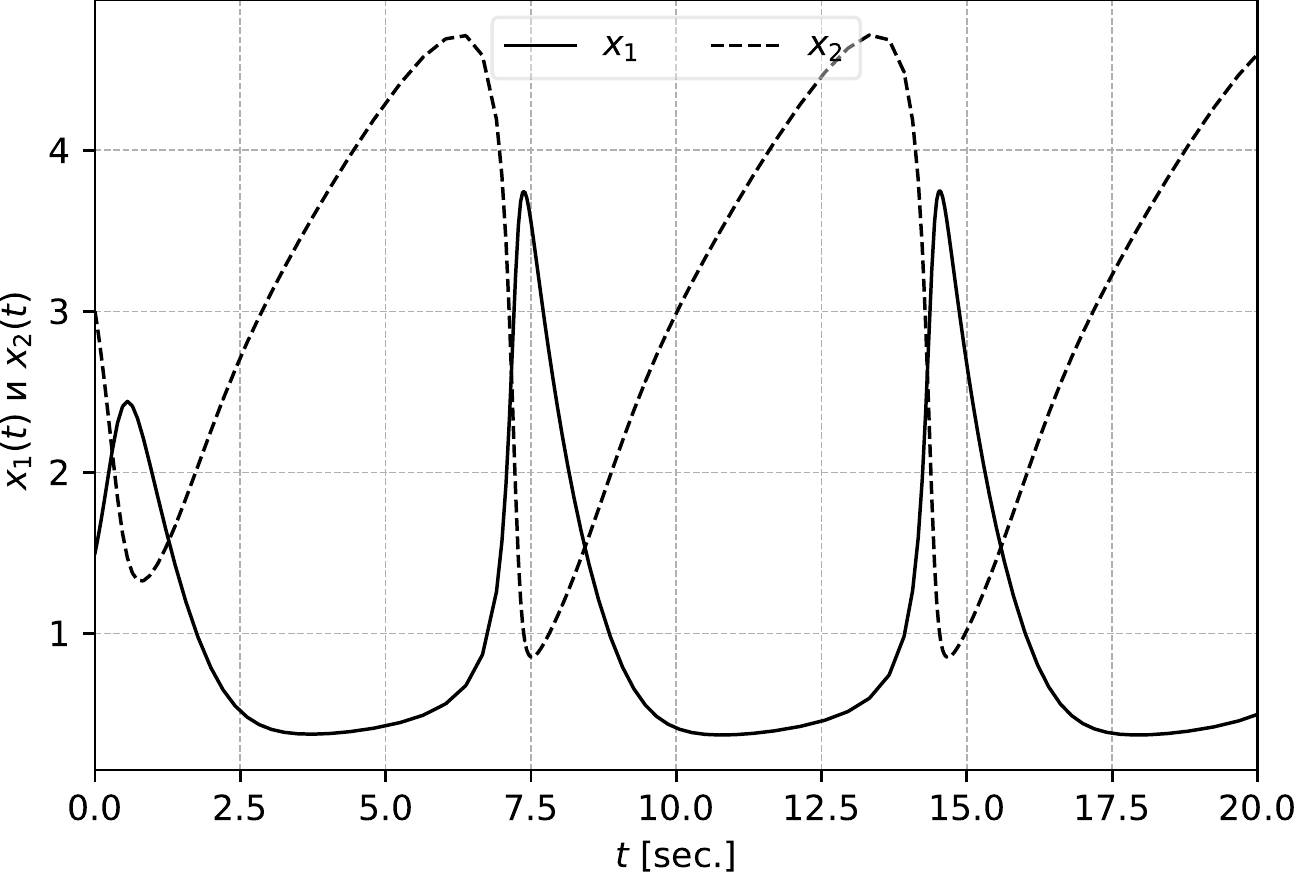}
\caption{Solution of system~\eqref{eq:01}}
  \label{fig:erk43b:01}
\end{figure}

\begin{figure}[!h]
  \centering
  \includegraphics[width=1.0\linewidth]{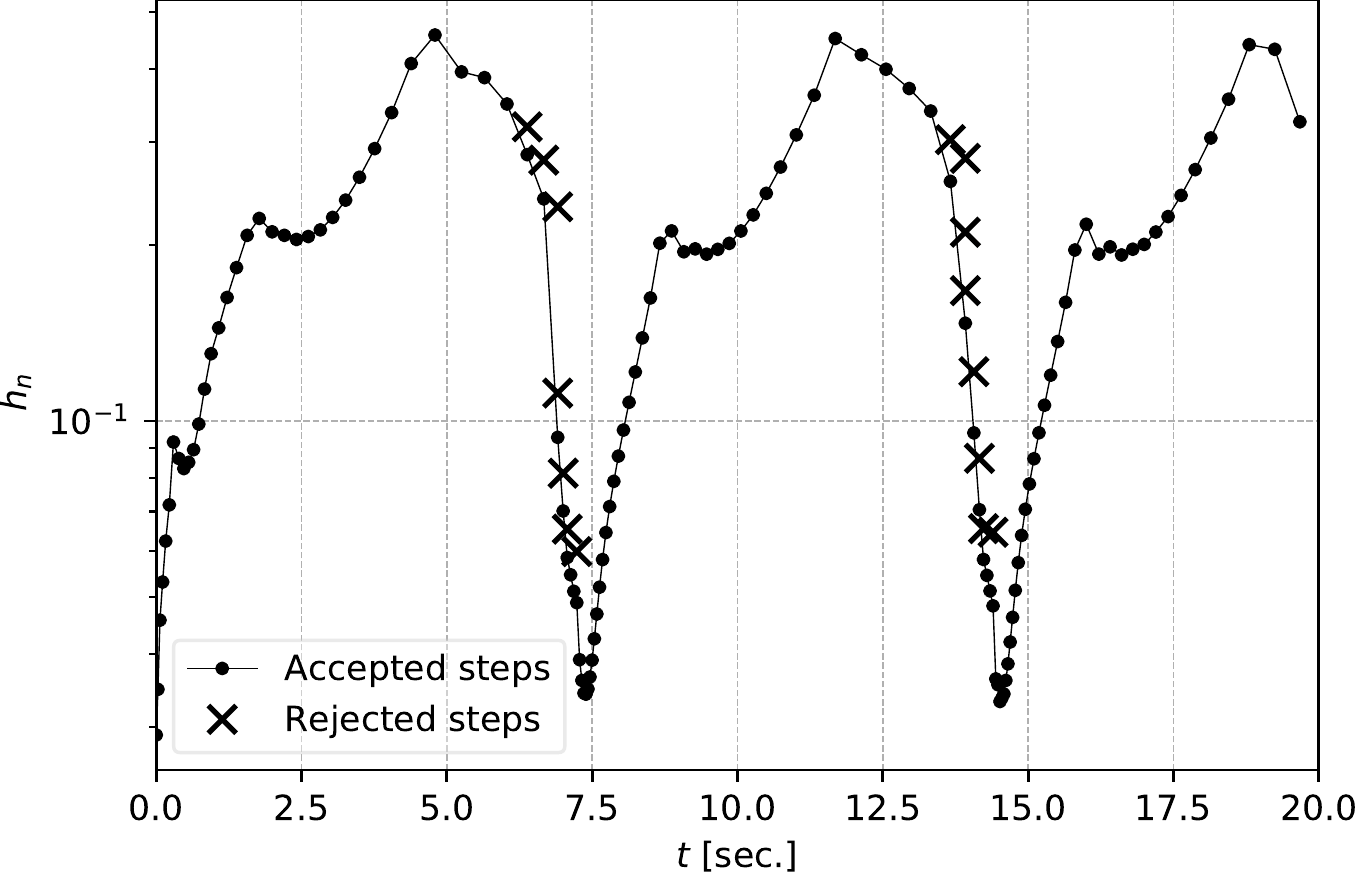}
\caption{Rejected and accepted steps for~\eqref{eq:01}}
  \label{fig:erk43b:02}
\end{figure}

Under the initial conditions $x_{1}(0) = 1.5$ and $x_{2}(0) = 3$, the Brusselator equation is
numerically integrated by the \verb|ERK43b| method on a segment $0 \leqslant t \leqslant 20$ with the
absolute and relative tolerances $A_{tol} = R_{tol} = 10^{-4}$. The corresponding graphs, which are
identical to those from~\cite[p.~170, Fig.~4.1]{L_Hairer_1}, are constructed. When
comparing the results from~\cite{L_Hairer_1} (see Fig.~\ref{fig:erk43b:hairer}) and our results (see
Figs.~\ref{fig:erk43b:01} and~\ref{fig:erk43b:02}), it can be seen that our method performs almost
identically to that from~\cite{L_Hairer_1}; however, our method selects step size
more accurately.

Another system of equations often used to test numerical schemes is a
special case of the restricted three-body problem: the problem of
computing Arenstorf orbits (named after an American physicist who
computed a stable orbit of a small body between the Moon and the
Earth).

The restricted three-body problem considers the motion of a small body
in the gravitational fields of a medium body and large body (the Moon
and the Earth). The mass of the small body is considered zero, while
the masses of the medium and large bodies are $\mu_1$
and $\mu_2$, respectively. The
orbits of the bodies are assumed to be on the same plane.

In the Arenstorf problem, different stable orbits can be obtained
depending on initial values~\cite{L_Hairer_1}. In dimensionless synodic coordinates,
three groups of initial values, which provide three different orbits,
have the following form:
\begin{equation}
  \label{eq:arenstorf01:initial}
  \begin{gathered}
    p_y = -1.00758510637908238, \quad p_x = 0.0,
    \\
    q^x = 0.994, \quad q^y =0.0;
  \end{gathered}
\end{equation}

\begin{equation}
  \label{eq:arenstorf02:initial}
  \begin{gathered}
    p_y = -1.03773262955733680, \quad p_x = 0.0,
    \\
    q^x = 0.994, \quad q^y = 0.0;
  \end{gathered}
\end{equation}

\begin{equation}
  \label{eq:arenstorf03:initial}
  \begin{gathered}
    p_y =  0.15064248999999985, \quad p_x = 0.0,
    \\
    q^x = 1.2, \quad q^y = 0.0.
  \end{gathered}
\end{equation}

Here, $p_x$ and $p_y$ are generalized momenta of the system, while $q^x$ and $q^y$ are
generalized coordinates of the system.

It should be noted that the initial values are purposefully set with
high accuracy because the small body very closely approaches the
medium body and even a small computational error can lead to an
incorrect physical interpretation (fall of the small body on the
medium body). This makes the Arenstorf problem well suitable for
testing implementations of numerical schemes.

The Hamilton function in synodic coordinates is written as follows:
\begin{equation*}
  H(p_x, p_y, q^x, q^y) = \frac12 (p_x^2 + p_y^2) + p_x q^y - p_y q^x - F(q^x, q^y),
\end{equation*}
where
\begin{gather*}
  F(q^x, q^y) = \frac{\mu_1}{r_1} + \frac{\mu_2}{r_2},
  \quad
  \mu_1 + \mu_2 = 1,
  \\
  r_1 = \sqrt{(q^x-\mu_2)^2 + (q^y)^2},
  \\
  r_2 = \sqrt{(q^x + \mu_1)^2 + (q^y)^2};
\end{gather*}
\begin{gather*}
  \pdv{F}{q^x} = - \frac{\mu_1 (q^x - \mu_2)}{r^3_1} - \frac{\mu_2
    (q^x + \mu_1)}{r^3_2},
  \\
  \pdv{F}{q^y} = - \frac{\mu_1 q^y}{r^3_1} - \frac{\mu_2 q^y}{r^3_2}.
\end{gather*}

The canonical equations to which the numerical scheme is applied have
the following form:
\begin{equation*}
  \left\{
    \begin{aligned}
      & \dot{p}_x = +p_y + \pdv{F}{q^x}, \dot{q}^x = p_x + q^y,\\
      & \dot{p}_y = -p_x + \pdv{F}{q^y}, \dot{q}^y = p_y - q^x.
    \end{aligned}
  \right.
\end{equation*}

\begin{figure}[!h]
  \begin{minipage}{\linewidth}
    \centering
    \includegraphics[width=0.7\linewidth]{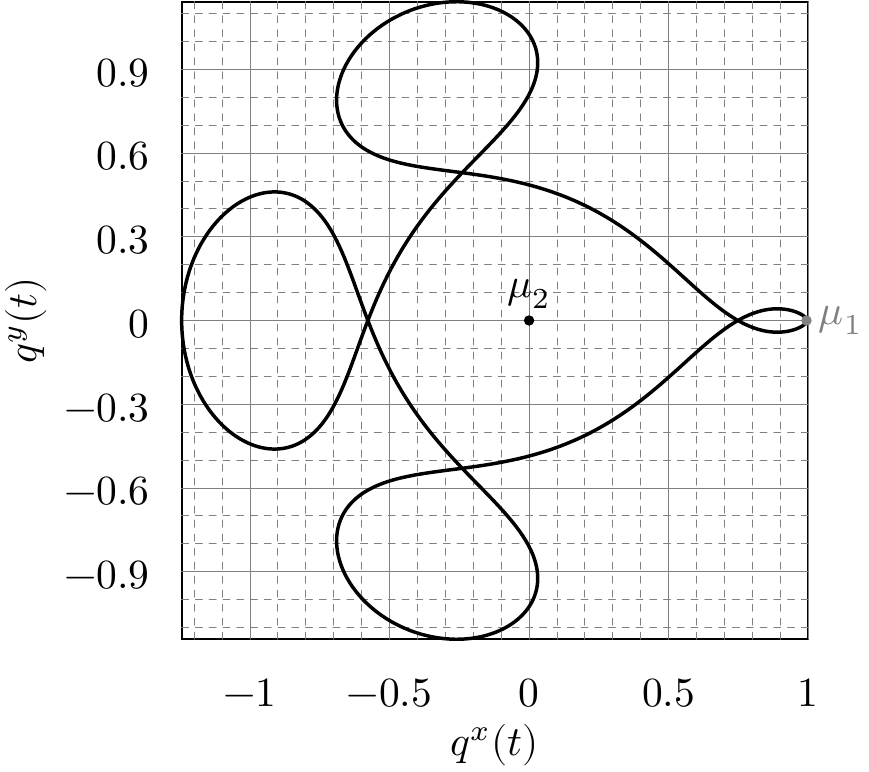}
    \caption{Orbit for the first group of initial values~\eqref{eq:arenstorf01:initial}}
    \label{fig:arenstorf01}
  \end{minipage}
  \\
  \begin{minipage}{\linewidth}
    \centering
    \includegraphics[width=0.7\linewidth]{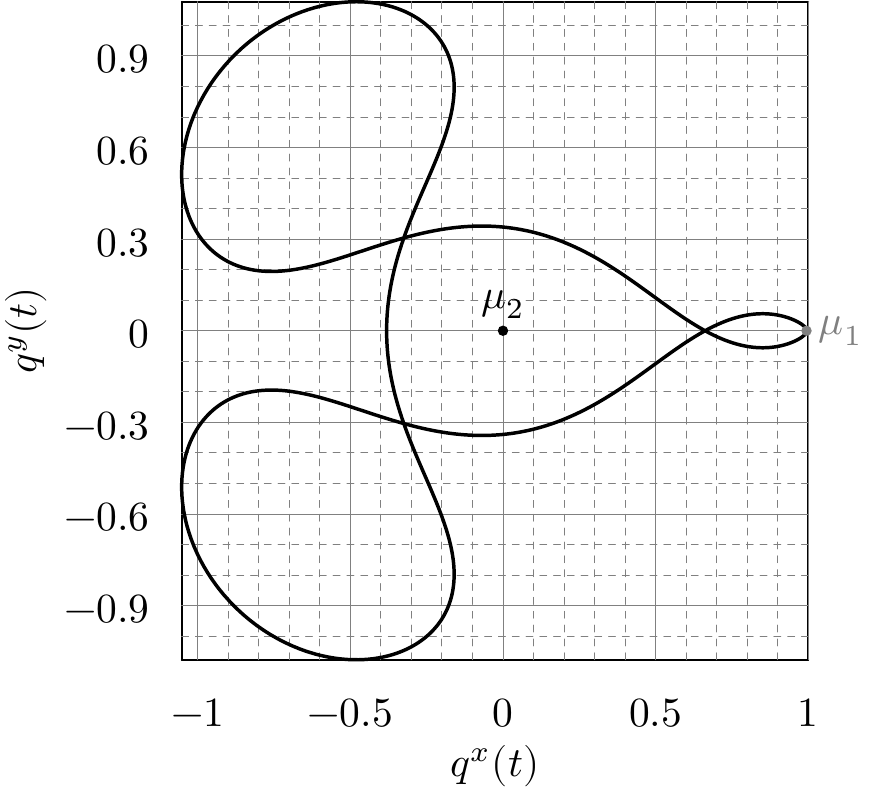}
\caption{Orbit for the second group of initial values~\eqref{eq:arenstorf02:initial}}
    \label{fig:arenstorf02}
  \end{minipage}
  \\
  \begin{minipage}{\linewidth}
    \centering
    \includegraphics[width=0.7\linewidth]{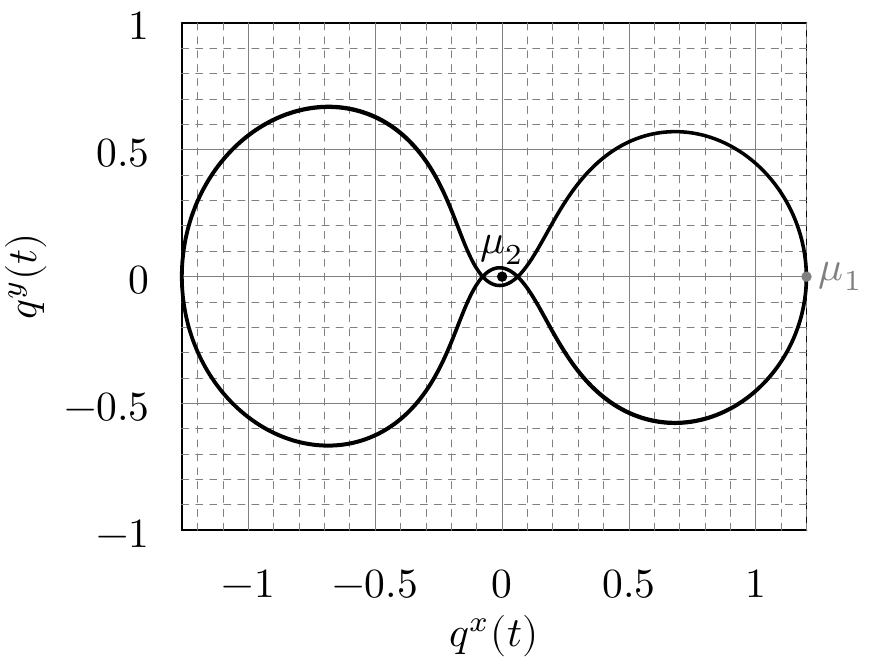}
\caption{Orbit for the third group of initial values~\eqref{eq:arenstorf03:initial}}
    \label{fig:arenstorf03}
  \end{minipage}
\end{figure}

Numerical solution is carried out for $0 \leqslant t \leqslant 17.065216560157962558$ with the absolute and relative
tolerances $A_{tol} = 10^{-17}$ and $R_{tol} = 0$ and average mass $\mu_1 = 0.012277471$. As a result, we find Arenstorf
orbits (see Figs.~\ref{fig:arenstorf01},
\ref{fig:arenstorf02}, and \ref{fig:arenstorf03}) for initial values~\eqref{eq:arenstorf01:initial},
\eqref{eq:arenstorf02:initial}, and~\eqref{eq:arenstorf03:initial}, respectively.

To estimate the error, the system is solved numerically for a time
interval equal to one period. At the end point of the interval, the
small body must return to the start point, i.e., $(q^x(0),q^y(0)) = (q^x(T),q^y(T))$, which is why the
error $\| \mathbf{q}_{0} - \mathbf{q}_{n} \|$ corresponds to the global error of the method. The errors for
the first group of initial values are shown in Table~\ref{tab:errors}.

\begin{table}[!h]
  \caption{Global errors of the numerical methods for one orbital revolution for the first group of initial values}
  \label{tab:errors}
  \begin{center}
    \begin{tabular}{ll}
      \hline
      \textbf{Method} & \textbf{Error}
\\                        
      \hline
      \verb|DPRK546S| & $1.84566 \cdot 10^{-12}$\\
      \verb|DPRK547S| & $2.93634 \cdot 10^{-12}$\\
      \verb|DPRK658M| & $4.22771 \cdot 10^{-12}$\\
      \verb|Fehlberg45| & $7.42775 \cdot 10^{-12}$\\
      \verb|DOPRI5| & $1.95463 \cdot 10^{-13}$\\
      \verb|DVERK65| & $1.67304 \cdot 10^{-12}$\\
      \verb|Fehlberg78B| & $5.45348 \cdot 10^{-12}$\\
      \verb|DOPRI8| & $1.06343 \cdot 10^{-11}$\\
      \hline
    \end{tabular}
  \end{center}
\end{table}
\section{Conclusions}  

Thus, the main contribution of this work is as follows. 
\begin{enumerate}
\item The language of the Jinja2 template engine has been used as a
  specialized symbolic computation system. As compared to a universal
  symbolic computation system, this approach has allowed us to make
  the software package simpler and more compact, as well as improve
  its portability.
  
\item A set of Python scripts has been created using the Jinja2
  language; these scripts generate Julia codes that implement
  numerical schemes of Runge--Kutta methods without step-size control,
  embedded methods with step-size control, and Rosenbrock methods with
  step-size control (all the routines are available at
  \url{https://bitbucket.org/mngev/rungekutta_generator}).
  
\item The generated functions have been tested by solving several
  typical problems (in this paper, we have considered in detail the
  Arenstorf problem as it depends more heavily on the accuracy of a
  numerical scheme); the module for the Julia language is available at
  \url{https://bitbucket.org/mngev/rangekutta-autogen}.
  
\end{enumerate}

\begin{acknowledgments}

The publication has been prepared with the support of the ``RUDN University Program 5-100''
and funded by Russian Foundation for Basic Research (RFBR) according to the research project
No~19-01-00645.

\end{acknowledgments}

 \bibliographystyle{elsarticle-num}

\bibliography{bib/symbolic-rk/cite}

\end{document}

\graphicspath{{image/symbolic-rk/ru/}{image/symbolic-rk/}{image/}}

\title{Использование шаблонизатора как инструментария компьютерной алгебры}

\author{М. Н. Геворкян}
\email{gevorkyan-mn@rudn.ru}
\affiliation{Кафедра прикладной информатики и теории вероятностей\\
  Российский университет дружбы народов\\
  ул. Миклухо-Маклая, д. 6, Москва, 117198, Россия}

\author{А. В. Королькова}
\email{korolkova-av@rudn.ru}
\affiliation{Кафедра прикладной информатики и теории вероятностей\\
  Российский университет дружбы народов\\
  ул. Миклухо-Маклая, д. 6, Москва, 117198, Россия}

\author{Д. С. Кулябов}
\email{kulyabov-ds@rudn.ru}
\affiliation{Кафедра прикладной информатики и теории вероятностей\\
  Российский университет дружбы народов\\
  ул. Миклухо-Маклая, д. 6, Москва, 117198, Россия}
\affiliation{Лаборатория информационных технологий\\
  Объединённый институт ядерных исследований\\
  ул. Жолио-Кюри, д. 6, Дубна, Московская область, 141980, Россия}

\begin{abstract}
  В исследовательских задачах, требующих применения численных методов
  решения систем обыкновенных дифференциальных уравнений, часто
  возникает необходимость выбора наиболее эффективного и оптимального
  для конкретной задачи численного метода. В частности, для решения
  задачи Коши, сформулированной для системы обыкновенных
  дифференциальных уравнений, применяются методы Рунге--Кутты (явные
  или неявные, с управлением шагом сетки или без и т.д.).  При этом
  приходится перебирать множество реализаций численного метода,
  подбирать коэффициенты или другие параметры численной схемы.  В
  данной статье предложено описание разработанной авторами библиотеки
  и скриптов автоматизации генерации функций программного кода на
  языке Julia для набора численных схем методов Рунге--Кутты. При этом
  для символьных манипуляций использовано программное средство
  подстановки по шаблону. Предлагаемый подход к автоматизации
  генерации программного кода позволяет вносить изменения не в каждую
  подлежащую сравнению функцию по отдельности, а использовать для
  редактирования единый шаблон, что с одной стороны дает
  универсальность в реализации численной схемы, а с другой позволяет
  свести к минимуму число ошибок в процессе внесения изменений в
  сравниваемые реализации численного метода.  Рассмотрены методы
  Рунге--Кутты без управления шагом, вложенные методы с управлением
  шагом и методы Розенброка также с управлением шагом. Полученные
  автоматически с помощью разработанной библиотеки программные коды
  численных схем протестированы при численном решении нескольких
  известных задач.
  % \textbf{Предпосылки}
  % В ходе учебной и исследовательской деятельности в
  % области численного решения систем обыкновенных дифференциальных
  % уравнений, возникла необходимость в переключении между многими
  % численными методами для изучения их эффективности и выбора
  % оптимального для той или иной задачи.
  % \textbf{Цель}
  % Создание библиотеки вложенных методов Рунге--Кутты с
  % управлением шагом на на языке программирования Julia, с возможностью
  % использования любой таблицы коэффициентов.
  % \textbf{Методы}
  % Автоматическая генерация программного кода.
  % \textbf{Результаты}
  % Созданы сценарии на языке Python с использованием
  % шаблонизатора Jinja2 для генерирования кода
  % вложенных методов Рунге--Кутты на языке Julia. Создан модуль для
  % языка Julia на основе сгенерированных методов.
\end{abstract}

  \keywords{}

\maketitle

\section{Введение}
\label{sec:intro}

Методы Рунге--Кутты являются основными численными методами для решения
нежестких систем обыкновенных дифференциальных уравнений. Известны
численные схемы высоких порядков (10 и выше) с управлением шагом и
плотной выдачей. Для многих языков программирования разработаны
библиотеки, реализующие наиболее эффективные из построенных численных
схем. Наиболее отлаженные коды были созданы для языка Fortran 77
Э.~Хайрером (Ernst Hairer) и его коллегами, они описаны в
книгах~\cite{L_Hairer_1,L_Hairer_2} и доступны для скачивания на
сайте~\cite{L_FortranCodes}. Эти программы реализуют методы
\texttt{DOPRI5}~\cite{L_DORPRI1980} и
\texttt{DOPRI853}~\cite{L_DORPRI1981}, которые обеспечивают управление
шагом и плотную выдачу. Второй из этих методов также поддерживает
переключение порядка метода между $8$ и $5$ для обеспечения большего
быстродействия при решении гладких задач.

Подпрограммы \texttt{DOPRI5} и \texttt{DOPRI853} очень хорошо
оптимизированы и доступны для использования во многих библиотеках и
математических пакетах, таких как Matlab~\cite{L_Matlab},
Octave~\cite{L_Octave}, SciPy~\cite{L_Scipy}, SciLab~\cite{L_Scilab},
Boost~C++~\cite{L_Boost} и т.д.
Однако часто для учебных и исследовательских целей возникает
необходимость перебора нескольких методов Рунге--Кутты для выбора
оптимального для конкретной
задачи~\cite{kulyabov:2016:rk-stochastic,kulyabov:2017:ceur:vol2064:runge-kutta,kulyabov:2018:dccn:rk-stochastic}. % Кроме
Кроме того, эффективность той или иной численной схемы также
существенно зависит от решаемой задачи.

Чтобы можно было достоверно сравнить разные численные схемы между
собой, необходимо иметь универсальный код, который реализует вложенные
методы для любого набора коэффициентов.  При этом все  реализации
должны быть унифицированы и отличаться лишь набором коэффициентов.

С целью решения поставленной задачи нами разработана библиотека для
языка программирования
Julia~\cite{L_JuliaLang_1,L_JuliaLang_2,L_JuliaLang_3}, описание
которой и приводится в данной статье. % Данный язык сравнительно новый,
Библиотека состоит из двух частей: вычислительной части, реализующей
алгоритмы типа Рунге--Кутты, и конструирующей части, создающей с
помощью символьных вычислений конкретный вариант алгоритма
Рунге--Кутты.  Отличительной особенностью конструирующей части
является автоматическая генерация кода функций из нескольких готовых
шаблонов.  Первоначально рассматривался вариант реализации
конструирующей части с помощью универсальной системе компьютерной
алгебры, например, с использованием
SymPy~\cite{kulyabov:2020:rng-cas-sympy::ru}. Однако в ходе реализации
выяснилось, что нам нужно от системы компьютерной алгебры в основном
операции сопоставления по образцу. В результате было принято решение
использовать более легковесное средство символьных
манипуляций. Фактически, для символьных манипуляций мы используем не
развитую систему компьютерной алгебры, а программное средство
подстановки по шаблону.  Такой подход хотя и несколько экзотичен, но
оправдан большей производительностью сгенерированных функций. Для
генерации кода использован язык Python~\cite{L_Python} и шаблонизатор
Jinja2~\cite{L_Jinja2}.

\subsection{Структура статьи}
\label{sec:structure}

В разделе~\ref{sec:rk} работы дается краткое описание вложенных
методов Рунге--Кутты и стратегии управления
шагом. Раздел~\ref{sec:program} посвящён описанию разработанной
библиотеки для языка программирования Julia. В разделе~\ref{sec:test}
приведены результаты тестирования сгенерированных автоматически с помощью
разработанной библиотек
функций для численного решения задач для систем обыкновенных дифференциальных
уравнений. В частности в этом разделе подробно рассмотрено решение
ограниченной задачи трех тел (задачи по вычислению орбит Аренсторфа).

\section{Методы Рунге--Кутты}
\label{sec:rk}

Методы Рунге--Кутты применяются к задаче Коши, сформулированной для
системы обыкновенных дифференциальных уравнений (ОДУ).  Теория методов
Рунге--Кутты хорошо известна и обстоятельно изложена в
книгах~\cite{L_Hairer_1, L_Hairer_2, L_Butcher}. Поэтому в данном
разделе мы приведем лишь основные формулы и кратко остановимся на
используемой в разработанной программной библиотеке стратегии
управления шагом.

\subsection{Постановка задачи}

Пусть даны две гладкие функции: неизвестная функция
$y^{\alpha}(t)\colon [t_0, T] \to \mathbb{R}^{N}$ и известная функция
$f^{\alpha}(t, y^{\beta}(t))\colon \mathbb{R}\times\mathbb{R}^N \to
\mathbb{R}^{N}$, где $[t_0, T] \in \mathbb{R}$,
$\alpha=1,\ldots,N$. Пусть также известно значение функции в начальный
момент времени $y^{\alpha}_0 = y^{\alpha}(t_0)$. Тогда задача Коши для
системы из $N$ обыкновенных дифференциальных уравнений формулируется
следующим образом:
\begin{equation}
  \label{Cauchy_ODE_0}
  \left\{
  \begin{aligned}
    &\dot{y}^{\alpha}(t) = f^{\alpha}(t, y^{\beta}(t)),\\
    & y^{\alpha}(t_0) = y^{\alpha}_0,\;\alpha,\beta = 1,\ldots,N.\\
  \end{aligned}
  \right.
\end{equation}

Систему~\eqref{Cauchy_ODE_0} можно записать покомпонентно:
\begin{equation}
  \label{Cauchy_ODE}
 \left\{
 \begin{aligned}
   &\dot{y}^1(t) = f^1(t, y^1(t),\ldots,y^N(t)),\\
  %  &\dot{y}^2(t) = f^2(t, y^1(t),\ldots,y^N(t)),\\
   &\vdots\\
   &\dot{y}^N(t) = f^N(t, y^1(t),\ldots,y^N(t)),\\
   & y^{\alpha}(t_0) = y^{\alpha}_0.\\ 
 \end{aligned}
 \right.
\end{equation}

Зададим на отрезке $[t_0, T]$ сетку из набора точек
$t_0 < t_1 < t_2 < \ldots < t_k < \ldots < t_N = T$ с шагом
сетки % Величина
$h_{k+1} = t_{k+1} - t_{k}$. % называется шагом сетки.
Каждой точке сетки по некоторому правилу, называемому \emph{численной
  схемой}, ставится в соответствие некоторая величина $y^{\alpha}_k$,
которая должна с необходимой точностью аппроксимировать решение
системы ОДУ в точках сетки, т.е.
$y^{\alpha}_k \approx y^{\alpha}(t_k)$. Для оценки погрешности
аппроксимации используют норму
$\| y^{\alpha}(t_{k}) - y^{\alpha}_{k} \|$.

\subsection{Явные вложенные методы Рунге--Кутты}

Погрешность метода Рунге--Кутты может оцениваться и другим
образом. Идея состоит в том, что кроме основного решения
$y^{\alpha}_{m}$ в точке рассматривается также вспомогательное решение
$\Hat{y}^{\alpha}_{m}$, выполненное с помощью метода Рунге--Кутты
смежного порядка. Разность этих решений может служить оценкой
локальной погрешности метода меньшего порядка.  Полученная оценка
локальной погрешности помогает выбрать переменный шаг интегрирования.
Методы, использующие данный способ оценки локальной погрешности,
называются \emph{вложенными} методами Рунге--Кутты~\cite{L_Hairer_1}.
При этом метод $y^{\alpha}_{m}$ называется  \emph{основным}, а метод 
$\Hat{y}^{\alpha}_{m}$~--- \emph{вложенным}.

Явный вложенный численный метод Рунге--Кутты для задачи
Коши~\eqref{Cauchy_ODE} задается следующими формулами:
\begin{equation*}
  \begin{aligned}
    & k^{1\alpha} = f^{\alpha}(t_{m}, y^{\beta}_{m}),
    \\
    &
    \begin{multlined}
      k^{i\alpha} = f^{\alpha}\big(t_{m} + c^{i}h_{m}, y^{\alpha}_{m}
      + h\sum^{i-1}_{j=1}a^{i}_{j}k^{j \beta}\big),
      \\
      i=2,\ldots,s,
    \end{multlined}
    \\
    &
    \begin{multlined}
      y^{\alpha}_{m+1}= y^{\alpha}_{m} + h_{m}(b_{1}k^{1\alpha} +
      b_{2}k^{2\alpha}
      + {} \\ {} +
      \ldots + b_{s-1}k^{s-1,\alpha} +
      b_{s}k^{s\alpha}),
    \end{multlined}
    \\
    &
    \begin{multlined}
      \hat{y}^{\alpha}_{m+1}= y^{\alpha}_{m} +
      h_{m}(\hat{b}_{1}k^{1\alpha} + \hat{b}_{2}k^{2\alpha}
      + {} \\ {} +
      \ldots +
      \hat{b}_{s-1}k^{s-1,\alpha} + \hat{b}_{s}k^{s\alpha}),
    \end{multlined}
  \end{aligned}
\end{equation*}
где $N$ --- число уравнений в системе ОДУ, $s$ --- число \emph{стадий}
численного метода. Латинские индексы $i,j,l=1,\ldots,s$ относятся к
численной схеме, греческие $\alpha,\beta$ --- к системе ОДУ, а
индексом $m$ обозначен номер шага.

Порядок аппроксимации решения сеточными функциям $y^{\alpha}_{m}$ и
$\hat{y}^{\alpha}_{m}$ различен, что позволяет на каждом шаге
вычислять погрешность аппроксимации, поэтому при указании порядка
вложенного метода используют обозначение $p(\Hat{p})$, где $p$ и
$\Hat{p}$ --- порядки основного и вложенного метода Рунге--Кутты
соответственно.

Величины $c^{i},a^{i}_{j},b_{j},\hat{b}_{j}$ полностью определяют
численную схему и называются \emph{коэффициентами метода},  их принято
группировать в виде \emph{таблицы Бутчера}, названной так в честь
Дж.~Ч.~Бутчера (John~C.~Butcher):
\begin{equation*}
  \begin{array}{c|cccccc}
  0  & 0 & 0 &\ldots & 0 & 0\\
  c^2  & a^{2}_1 & 0 &\ldots & 0 & 0\\
  c^3  & a^{3}_1 &a^{3}_2 &\ldots &0 & 0\\
  \vdots  &\vdots  &\vdots &\ddots &\vdots &\vdots\\
  c^{s-1} & a^{s-1}_1 &a^{s-1}_2 &\ldots & 0 & 0\\
  c^{s}  & a^{s}_1 &a^{s}_2 &\ldots &a^{s}_{s-1} & 0\\ \hline
       &b_{1} &b_{2} &\ldots &b_{s-1} & b_{s}\\
       &\hat{b}_{1} &\hat{b}_{2} &\ldots &\hat{b}_{s-1} & \hat{b}_{s}
  \end{array}
\end{equation*}

Нахождение коэффициентов метода Рунге--Кутты для порядков выше $4$
является отдельной сложной задачей, которой мы здесь не будем касаться
и которая подробно раскрыта в книге~\cite{L_Hairer_1}. Заметим только,
что на коэффициенты $c^{i}$ и $a^{i}_{j}$ также принято налагать
дополнительные условия следующего вида:
$c^{i} = a^{i}_{1} + \ldots + a^{i}_{s}$.

Явный метод Рунге-Кутты наиболее прост для реализации в программном виде, так как
на каждом шаге метода необходимо последовательно вычислять элементы
двухмерного массива \texttt{k[1:s,1:N]} в цикле по первому измерению
\verb|i=1:s|. Элемент \texttt{k[1,1:N]} требует только значений
$t_{m}, y_{m}$, вычисленных на предыдущем шаге, элемент
\texttt{k[2,1:N]} требует \texttt{k[1,1:N]}, элемент \texttt{k[3,1:N]}
требует \texttt{k[1,1:N]} и \texttt{k[2,1:N]} и т.д. Заметим, что
индексация элементов массива в языке Julia начинается с $1$, также как
и в языке Fortran.

К настоящему времени найдены коэффициенты для множества вложенных
методов Рунге--Кутты вплоть до 10-го порядка точности. Наиболее часто
используются методы, найденные Дж.\,Р.~Дорманом в соавторстве с
П.\,Дж.~Принцем (J.\,R.~Dorman,
P.\,J.~Prince)~\cite{L_DORPRI1980,L_DORPRI1981,L_DORPRI1986},
Э.~Фельбергом
(E.~Fehlberg)~\cite{fehlberg:1969:klassische-runge-kutta-funfter,fehlberg:1970:klassische-runge-kutta-vierter}
и Дж.\,Р.~Кэшом в соавторстве с А.\,Х.~Карпом (J.\,R.~Cash,
A.\,H.~Karp)~\cite{L_CashKarp1990}.

\subsection{Стратегии управления длиной шага}

Существуют различные стратегии управления величиной шага $h$ в
зависимости от локальной погрешности. Выбор той или иной стратегии
чаще всего диктуется характером решаемой задачи. Здесь мы изложим
классический алгоритм, предложенный в книге
Э.~Хайрера~\cite{L_Hairer_2}, который основан на идеях из теории
управления и хорошо работает для большинства нежестких задач.

Пусть $y^{\alpha}_{m}$ и $\hat{y}^{\alpha}_{m}$ --- два численных
решения разного порядка аппроксимации, вычисленные на шаге $m$ 
работы алгоритма. Потребуем, чтобы
$|y^{\alpha}_{m} - \hat{y}^{\alpha}_{m}| \leqslant sc^{\alpha}$, где
величина желаемой локальной погрешности $sc^{\alpha}$ (scale)
вычисляется по формуле
\begin{equation*}
  sc^{\alpha} = A_{tol} + \max(|y^{\alpha}_{m}|, |\hat{y}^{\alpha}_{m}|)R_{tol},
\end{equation*}
где $A_{tol}$ и $R_{tol}$ --- желаемая величина абсолютной и
относительной погрешности соответственно. Допустимая ошибка на шаге $m$ находится
как среднее квадратичное:
\begin{equation*}
  E_{m} = \sqrt{\dfrac{1}{N} \sum\limits_{\alpha = 1}^{N}{\left(\dfrac{y^{\alpha}_{m} - y^{\alpha}_{m}}{sc^{\alpha}}\right)^{2}}}.
\end{equation*}

Новая величина шага вычисляется по формуле
\begin{equation*}
  h_{m+1} = h_{m}/\max\left(f_{\min}, \min\left(f_{\max}, E^{a}_{m}E^{-b}_{m-1}/f_{s}\right)\right),
\end{equation*}
где $E_{m-1}$ --- величина ошибки, вычисленная на предыдущем
шаге. Фактор $f_{s}$ обычно равен $0.9$ или $0.8$ и предназначен для
предотвращения слишком резкого увеличения величины шага. Показатели
степени $a$ и $b$ подбираются исходя из задачи. В
книге~\cite{L_Hairer_2} рекомендуются следующие универсальные значения
$a = 0.7/p - 0.75b$ и $b = 0.4/p$, где $p$ --- порядок аппроксимации
используемого метода. Факторы $f_{\min}$ и $f_{\max}$ позволяют
ограничить границы изменения размера шага и также зависят от решаемой
задачи. На практике хорошо зарекомендовали себя значения
$f_{\min} = 0.1$ и $f_{\max} = 5.0$.

В случае $E_{m} < 1$ вычисленные значения $y^{\alpha}_{m}$ считаются
удовлетворительными и метод переходит к следующему шагу
вычислений. Если же $E_{m} > 1$, то вычисления считаются
неудовлетворительными и текущий шаг повторяется с новым значением
величины $h_{m}$
\begin{equation*}
  h_{m} \leftarrow h_{m} / \min\left(f_{\max}, E^{a}_{m}/f_{s}\right).
\end{equation*}

В качестве начального значения $y^{\alpha}_{m}$ на следующем шаге
$m+1$ можно использовать как непосредственно $y^{\alpha}_{m}$, так и
значение вложенного метода $\hat{y}^{\alpha}_{m}$.

\section{Описание программы}
\label{sec:program}

\subsection{Мотивация использования генерации кода}

Универсальная реализация явного вложенного метода Рунге--Кутты
подразумевает, что программе в качестве параметров должны передаваться
массивы коэффициентов из таблиц Бутчера, которые она будет
использовать для проведения расчетов. Очевидным способом хранения
данных коэффициентов является использование массивов. Однако матрица
$a^{i}_{j}$ является нижне-диагональной и хранение ее в виде массива
$s\times s$ приводит к тому, что более чем половина выделенной для
массива памяти тратится на хранение нулей. Коэффициенты $c^{i}$,
$b_{j}$ и $\hat{b}_{j}$ также часто имеют в своем составе нули,
хранение которых неразумно.

В связи с этим большинство кодов, реализующих явные вложенные методы
Рунге--Кутты, используют для хранения коэффициентов набор именованных
констант, а не массивы. Это вызвано также и тем, что операции со
скалярными величинами в большинстве языков проводятся быстрее, чем
операции с массивами.

В языке Julia стандартные массивы являются динамическими, поэтому
накладные расходы на хранение двумерного массива $s\times s$ больше,
чем на расходы на хранение $s \times s$ именованных констант.

При сохранении требования универсальности создаваемого кода и вместе с
тем желание увеличить скорости вычислений и уменьшить расход памяти,
привели нас к решению использовать автоматическую генерацию кода по
одному шаблону для каждого отдельного метода.

Кроме выигрыша в производительности, автоматическая генерации кода
позволяет добавлять или изменять не каждую функцию по отдельности, а
все функции совокупно путём редактирования одного лишь шаблона. Это
позволяет как уменьшить количество ошибок, так и генерировать
различные варианты функций для разных целей.

\subsection{Описание реализации генератора методов Рунге--Кутты}

В качестве языка для генерации кода нами был выбран язык Python, так
как он поддерживает разнообразные средства для работы с текстом. Также
в стандартной библиотеке Python присутствует тип данных
\texttt{Fraction}, который позволяет задать коэффициенты Рунге--Кутты
в виде рациональных дробей, а потом уже преобразовывать их в
вещественный вид с нужным порядком точности.

Кроме самого языка Python была использована библиотека для обработки
шаблонов Jinja2~\cite{L_Jinja2}. Данный шаблонизатор (template engine)
разрабатывался изначально для генерации HTML-страниц, однако он
обладает очень гибким синтаксисом и может использоваться как
универсальное средство для генерации текстовых файлов любого вида, в
том числе и исходных кодов на любых языках программирования. Кроме
Jinja2 мы использовали библиотеку \texttt{numpy}~\cite{L_Scipy} для
работы с массивами коэффициентов.

Шаблоны для генерации функций хранятся в файлах \texttt{rk\_method.jl}
и \texttt{erk\_method.jl}, которые представляют собой исходный код на
языке Julia с инструкциями Jinja2. Шаблонизатор позволяет поместить
всю логику генерации кода в шаблон и передавать извне только данные о
методе.

Информация о численных схемах хранится в виде JSON-файла, где каждый
метод представляет собой JSON-объект примерно следующего вида:
\begin{minted}{json}
 { "name": "Название (одним словом)",
  "description": "Краткое описание",
  "stage": s,
  "order": p,
  "extrapolation_order": p_hat,
  "a": [["0", "0", "0"], ["1", "0", "0"], ["1/4", "1/4", "0"]],
  "b": ["1/2", "1/2", "0"],
  "b_hat": ["1/6", "1/6", "2/3"],
  "c": ["0", "1", "1/2"]}
\end{minted}

Название JSON-объекта используется в качестве имени сгенерированной в
последствии функции. Массивы \verb|a|, \verb|b|, \verb|b_hat|,
\verb|c| могут быть как числового, так и строкового типа. %В случае
Если коэффициенты метода заданы в виде рациональных дробей, то можно
указать их в виде \verb|"m/n"|, затем % . В дальнейшем
они будут преобразованы
в объект \verb|Fraction| стандартной библиотеки Python, а в теле
сгенерированной функции будут представлены десятичной дробью двойной
точности с $17$ значащими знаками.

Список объектов указанного выше вида последовательно обрабатывается скриптами
\texttt{erk\_generator.py} и \texttt{rk\_generator.py}. Скрипты,
используя переданную им информацию о методах, генерируют для каждого
метода код нескольких функций языка Julia.

Скрипты позволяют сгенерировать функции для $16$ вложенных методов
Рунге--Кутты. Методы низкого порядка взяты из книги
Хайрера~\cite{L_Hairer_1}. Из методов пятого порядка и выше заданы
коэффициенты из работ~\cite{L_DORPRI1980,L_DORPRI1981,L_DORPRI1986}
и~\cite{fehlberg:1969:klassische-runge-kutta-funfter,fehlberg:1970:klassische-runge-kutta-vierter}. Метод
из работы~\cite{L_CashKarp1990} реализован отдельно, так как
предусматривает переключение порядка точности и сложный алгоритм
управления шагом. Для добавлении своих методов следует задать
JSON-объект указанного выше формата и запустить скрипты вновь. Для
добавленных методов будут сгенерированы соответствующие функции.

\subsection{Описание сгенерированных функций}

Разработанная нами библиотека имеет стандартную для модулей языка
Julia структуру. Весь исходный код находится в каталоге \verb|src|. В
подкаталоге \verb|src/generated| располагаются файлы с автоматически
генерируемыми функциями, которые, в свою очередь, с помощью директивы
\verb|include| включаются в главный файл модуля
\verb|src/RungeKutta.jl|. Код, ответственный за управление шагом,
находится в файле \verb|StepControl.jl| и является общим для всех
методов.

Для каждого из вложенных методов Рунге--Кутты генерируются три функции,
две из которых имеют следующий вид:
\begin{minted}{julia}
  ERK(func, A_tol, R_tol, x_0, t_start, t_stop) -> (T, X)
  ERK(func, A_tol, R_tol, x_0, t_start, t_stop, last) -> (t, x)
\end{minted}

Функции имеют одинаковое имя, а список их аргументов различается лишь
последним аргументом. Благодаря поддержке языком Julia множественной
диспетчеризации (перегрузка функций), компилятор сам определяет в
каком случае какую из реализаций вызывать. 

В указанных выше функциях:
\begin{itemize}
\item \verb|func::Function| --- правая часть системы ОДУ
  $\dot{x}^{\alpha}(t) = f^{\alpha}(t, x^{\beta})$: \verb|func(t, x)|,
  где \verb|t::Float64|~--- время, \verb|x::Vector{Float64}|~--- значение
  функции $x(t)$; аргумент \verb|t| должен всегда присутствовать в
  вызове функции, даже если система автономная и явно от времени не
  зависит;
\item аргументы \verb|A_tol| и \verb|R_tol| должны иметь тип
  \verb|Float64| и обозначают абсолютную и относительную точности
  метода соответственно;
\item \verb|x_0::Vector{Float64}|~--- начальное значение функции
  $x^{\alpha}_{0} = x^{\alpha}(t_0)$;
\item \verb|t_start| и \verb|t_stop| имеют тип \verb|Float64| и
  обозначают начальную и конечную точки промежутка интегрирования;
\item если указан аргумент \verb|last::Bool|, то будут возвращены
  последняя точка промежутка интегрирования $t_n$ и вектор
  $x^{\alpha}_{n} \approx x^{\alpha}(t_n)$, в противном случае
  возвращаются массивы \verb|T::Vector{Float64}| и
  \verb|Matrix{Float64}|.
\end{itemize}

Для методов Рунге--Кутты без управления шагом генерируются аналогичные
функции, единственное отличие которых заключается в отсутствии
аргументов \verb|A_tol| и \verb|R_tol|. Вместо них следует передавать
аргумент \verb|h|, который задает размер шага, используемого для
вычислений.

Для реализации % изучения работы
алгоритма управления шагом генерируется еще одна
функция:
\begin{minted}{julia}
  ERK_info(func, A_tol, R_tol, x_0, t_start, t_stop) -> (accepted_t, accepted_h, rejected_t, rejected_h)
\end{minted}
Эта функция возвращает следующие массивы:
\begin{itemize}
\item \verb|accepted_t|: все точки сетки, в которых погрешность
  вычислений была признана удовлетворительной;
\item \verb|accepted_h|: все принятые размеры шагов;
\item \verb|rejected_t|: все точки сетки, в которых погрешность
  вычислений была признана \textbf{не}удовлетворительной;
\item \verb|rejected_h|: все отбракованные размеры шага;
\item \verb|errors|: значения локальной ошибки $E_{n}$.
\end{itemize}

Все возвращаемые значения имеют тип \verb|Vector{Float64}|.
  \section{Тестирование сгенерированных программных кодов численных  схем}
\label{sec:test}

Для тестирования созданных методов были использованы три системы
дифференциальных уравнений, рассмотренные в~\cite{L_Hairer_2}.  Первая
система --- уравнения ван дер Поля:
\begin{gather*}
  \dv[2]{x}{t} - \mu(1-x^2) \dv{x}{t} + x = 0,
  \\
  \left\{
  \begin{aligned}
      & \Dot{x}_{1}(t) = x_{2}(t),\\
      & \Dot{x}_{2}(t) = \mu (1-x^{2}_{1}(t)) x_{2}(t) - x_{1}(t),
  \end{aligned}
  \right.
\end{gather*}
заданные со следующими начальными значениями:
\begin{equation*}
  \mu = 1,
  \quad
  \mathbf{x} = (0, \sqrt{3})^T,
  \quad
  0 \leqslant t \leqslant 12,
\end{equation*}
где коэффициент $\mu$ характеризует нелинейность и силу затухания колебаний.

Вторая система --- уравнения твёрдого тела без внешних сил, или иначе,
уравнения Эйлера твердого тела:
\begin{equation*}
  \left\{
  \begin{aligned}
      & \dot{x}_{1}(t) = I_{1} x_{2}(t) x_{3}(t),\\
      & \dot{x}_{2}(t) = I_{2} x_{1}(t) x_{3}(t),\\
      & \dot{x}_{3}(t) = I_{3} x_{1}(t) x_{2}(t),\\
      & I_1 = -2,\;I_2 = 1.25,\; I_3 = -0.5,\\
      & \mathbf{x} = (0, 1, 1)^T,\; t = [0, 12].
  \end{aligned}
  \right.
\end{equation*}

Для проверки работы алгоритма управления шагом применялась численная
схема вложенного метода Рунге--Кутты с $p=4$ и $\Hat{p}=3$ к системе
уравнений брюсселятора:
\begin{equation}
  \label{eq:01}
  \left\{
  \begin{aligned}
      & \dot{x}_{1}(t) = 1 + x^{2}_{1}x_{2} - 4x_{1},\\
      & \dot{x}_{2}(t) = 3x_{1} - x^{2}_{1} x_{2}.
  \end{aligned}
  \right.
\end{equation}

\begin{figure}[!h]
  \centering
  \includegraphics[width=1.0\linewidth]{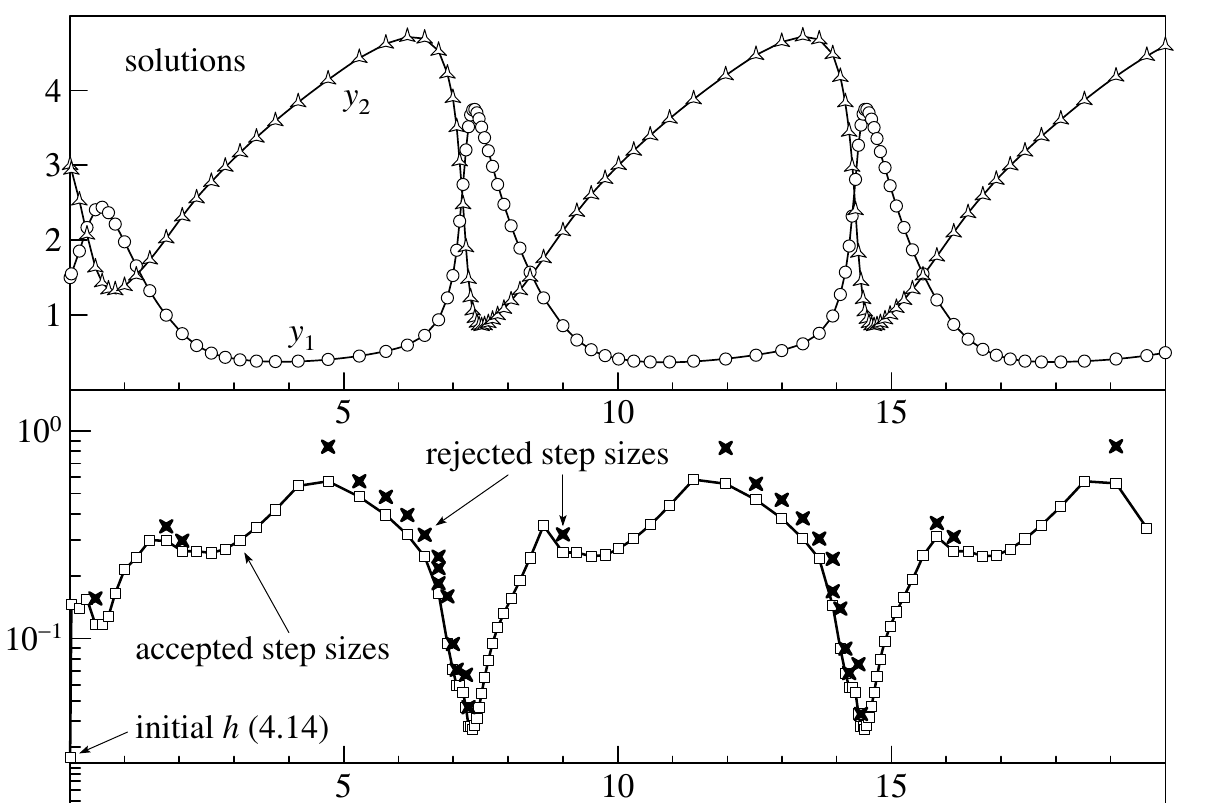}
  \caption{Решение системы~\eqref{eq:01} и отброшенные и принятые шаги (рисунок из~\cite[с.~170, рис.~4.1]{L_Hairer_1})} %
  \label{fig:erk43b:hairer}
\end{figure}

\begin{figure}[!h]
  \centering
  \includegraphics[width=1.0\linewidth]{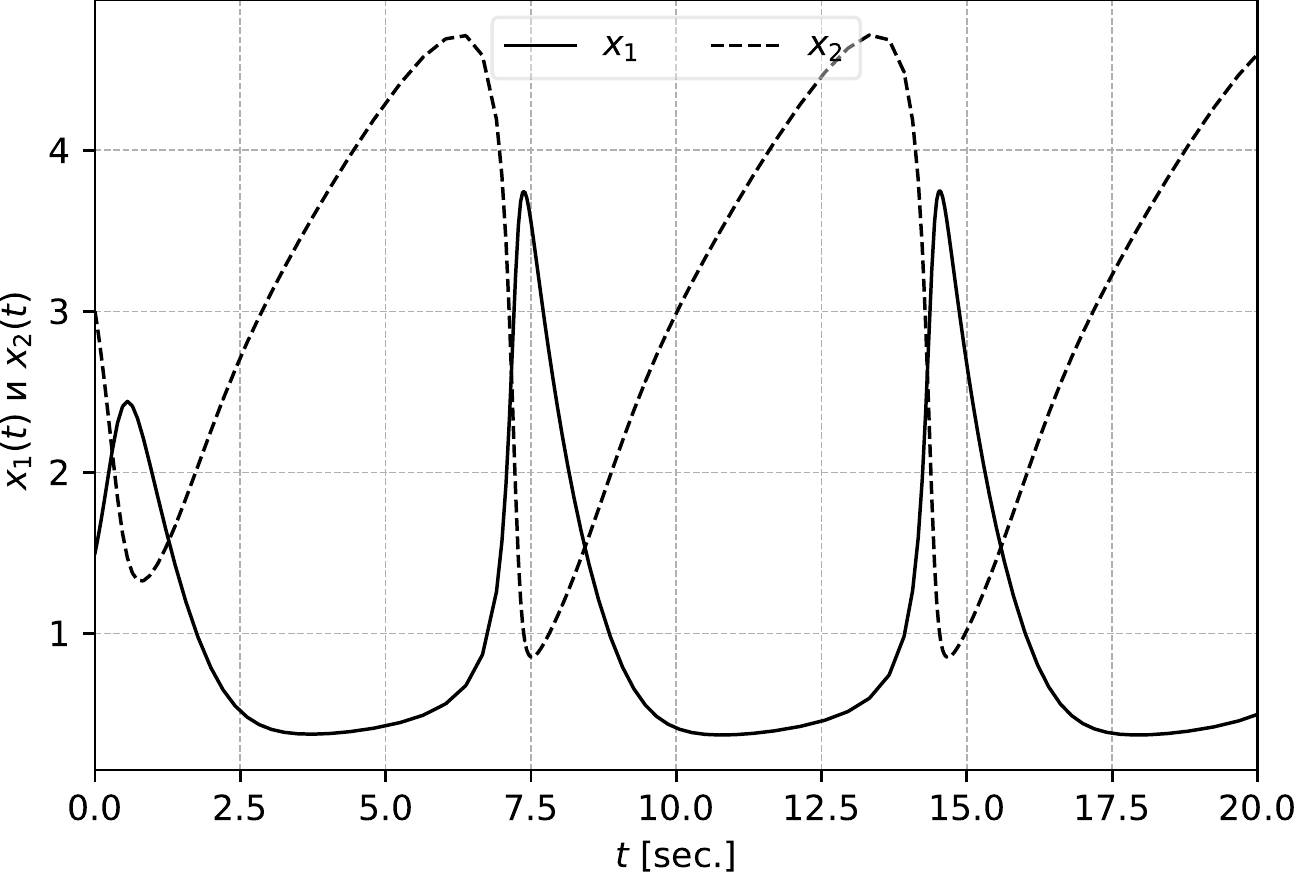}
  \caption{Решение системы~\eqref{eq:01}}
  \label{fig:erk43b:01}
\end{figure}

\begin{figure}[!h]
  \centering
  \includegraphics[width=1.0\linewidth]{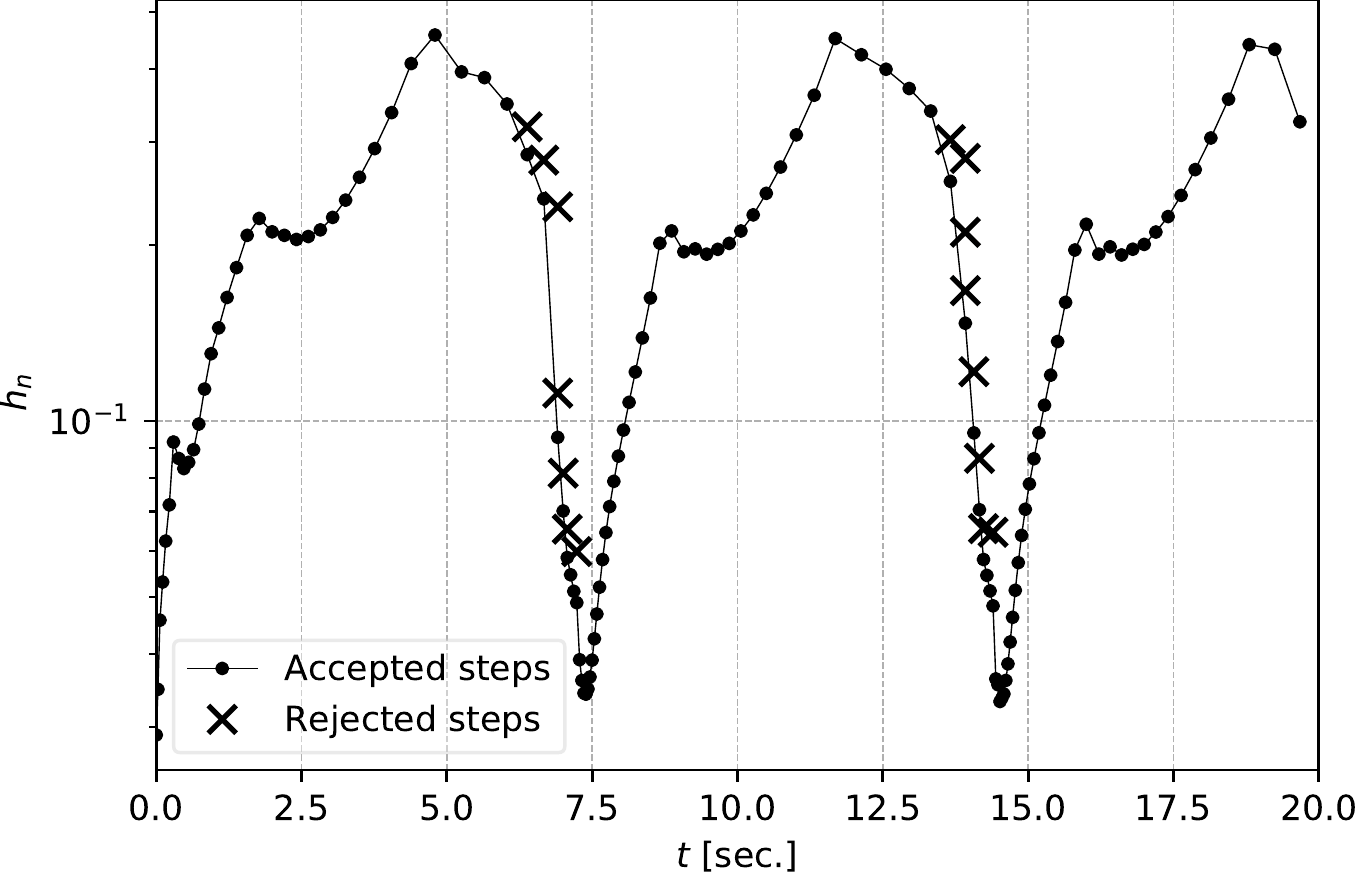}
  \caption{Отброшенные и принятые шаги~\eqref{eq:01}}
  \label{fig:erk43b:02}
\end{figure}

Для начальных условий $x_{1}(0) = 1.5$, $x_{2}(0) = 3$ уравнение
брюсселятора было численно проинтегрировано методом \verb|ERK43b| на
отрезке $0 \leqslant t \leqslant 20$ с абсолютной и относительной
погрешностями $A_{tol} = R_{tol} = 10^{-4}$. Были построены графики,
идентичные графикам из книги~\cite[с.~170, рис.~4.1]{L_Hairer_1}. При сравнении
результатов из книги~\cite{L_Hairer_1} (см. рис.~\ref{fig:erk43b:hairer}) с полученными
нами результатами
(см. рис.~\ref{fig:erk43b:01} и~\ref{fig:erk43b:02}) видно, что реализованный нами метод
работает практически также, как и использованный в~\cite{L_Hairer_1},
однако выбирает размер шага аккуратнее.

Ещё одна часто применяющаяся для тестирования численных схем система
уравнений --- частный случай ограниченной задачи трех тел, называемый
задачей по вычислению орбит Аренсторфа. Свое название данная задача
получила в честь Ричарда Ф. Аренсторфа~--- американского физика,
рассчитавшего стабильную орбиту малого тела между Луной и Землёй.

В ограниченной задаче трех тел рассматривается движение малого тела в
поле тяготения среднего и большого тел (Луна и Земля). Масса малого
тела считается равной нулю, массы среднего и большого тел~--- $\mu_1$ и
$\mu_2$ соответственно. Предполагается, что орбиты тел лежат в одной
плоскости.

В задаче Аренсторфа в зависимости от начальных значений можно получить
разные стабильные орбиты~\cite{L_Hairer_1}. В безразмерных синодических координатах три
группы начальных значений, дающих три разные орбиты, имеют следующий
вид:
\begin{equation}
  \label{eq:arenstorf01:initial}
  \begin{gathered}
    p_y = -1.00758510637908238, \quad p_x = 0.0,
    \\
    q^x = 0.994, \quad q^y =0.0;
  \end{gathered}
\end{equation}

\begin{equation}
  \label{eq:arenstorf02:initial}
  \begin{gathered}
    p_y = -1.03773262955733680, \quad p_x = 0.0,
    \\
    q^x = 0.994, \quad q^y = 0.0;
  \end{gathered}
\end{equation}

\begin{equation}
  \label{eq:arenstorf03:initial}
  \begin{gathered}
    p_y =  0.15064248999999985, \quad p_x = 0.0,
    \\
    q^x = 1.2, \quad q^y = 0.0.
  \end{gathered}
\end{equation}

Здесь $p_x$, $p_y$ --- обобщённые импульсы системы, $q^x$, $q^y$~---
обобщенные координаты системы.

Важно отметить, что начальные значения специально указаны с большой
точностью, так как малое тело крайне близко подходит к среднему телу и
даже небольшая погрешность вычислений может привести к ошибке, что
физически интерпретируется как падение малого тела на среднее. Как раз
из-за такой особенности данная задача хорошо подходит для проверки
реализации численных схем.

Функция Гамильтона в синодических координатах записывается следующим
образом:
\begin{equation*}
  H(p_x, p_y, q^x, q^y) = \frac12 (p_x^2 + p_y^2) + p_x q^y - p_y q^x - F(q^x, q^y),
\end{equation*}
где
\begin{gather*}
  F(q^x, q^y) = \frac{\mu_1}{r_1} + \frac{\mu_2}{r_2},
  \quad
  \mu_1 + \mu_2 = 1,
  \\
  r_1 = \sqrt{(q^x-\mu_2)^2 + (q^y)^2},
  \\
  r_2 = \sqrt{(q^x + \mu_1)^2 + (q^y)^2};
\end{gather*}
\begin{gather*}
  \pdv{F}{q^x} = - \frac{\mu_1 (q^x - \mu_2)}{r^3_1} - \frac{\mu_2
    (q^x + \mu_1)}{r^3_2},
  \\
  \pdv{F}{q^y} = - \frac{\mu_1 q^y}{r^3_1} - \frac{\mu_2 q^y}{r^3_2}.
\end{gather*}

Канонические уравнения, к которым и будет применяться численная схема,
имеют следующий вид:
\begin{equation*}
  \left\{
    \begin{aligned}
      & \dot{p}_x = +p_y + \pdv{F}{q^x}, \dot{q}^x = p_x + q^y,\\
      & \dot{p}_y = -p_x + \pdv{F}{q^y}, \dot{q}^y = p_y - q^x.
    \end{aligned}
  \right.
\end{equation*}

\begin{figure}[!h]
  \begin{minipage}{\linewidth}
    \centering
    \includegraphics[width=0.7\linewidth]{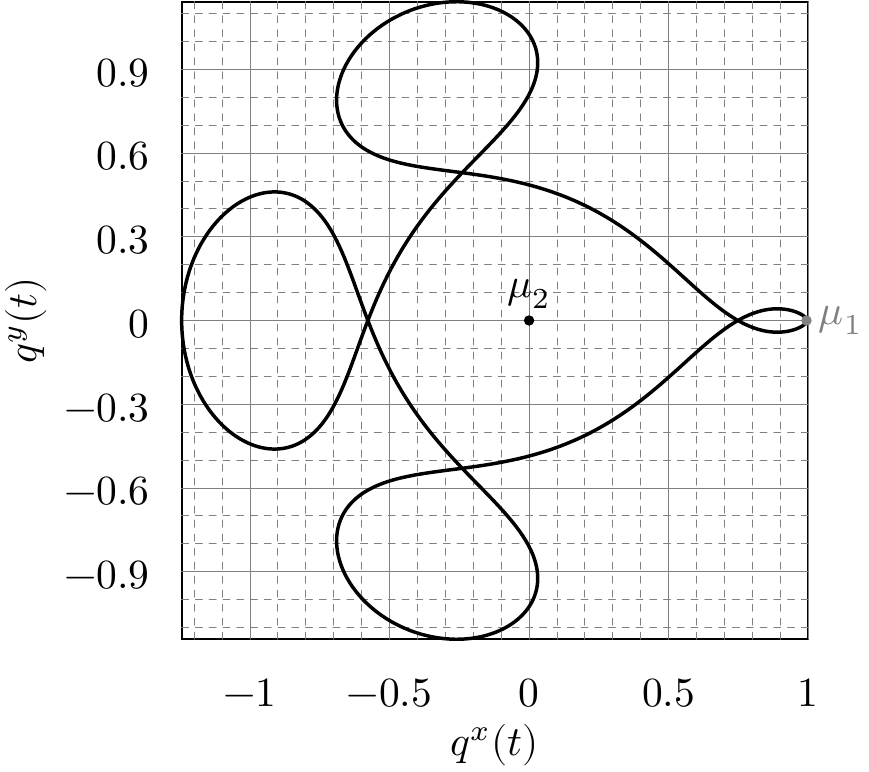}
    \caption{Орбита для первой группы начальных значений~\eqref{eq:arenstorf01:initial}}
    \label{fig:arenstorf01}
  \end{minipage}
  \\
  \begin{minipage}{\linewidth}
    \centering
    \includegraphics[width=0.7\linewidth]{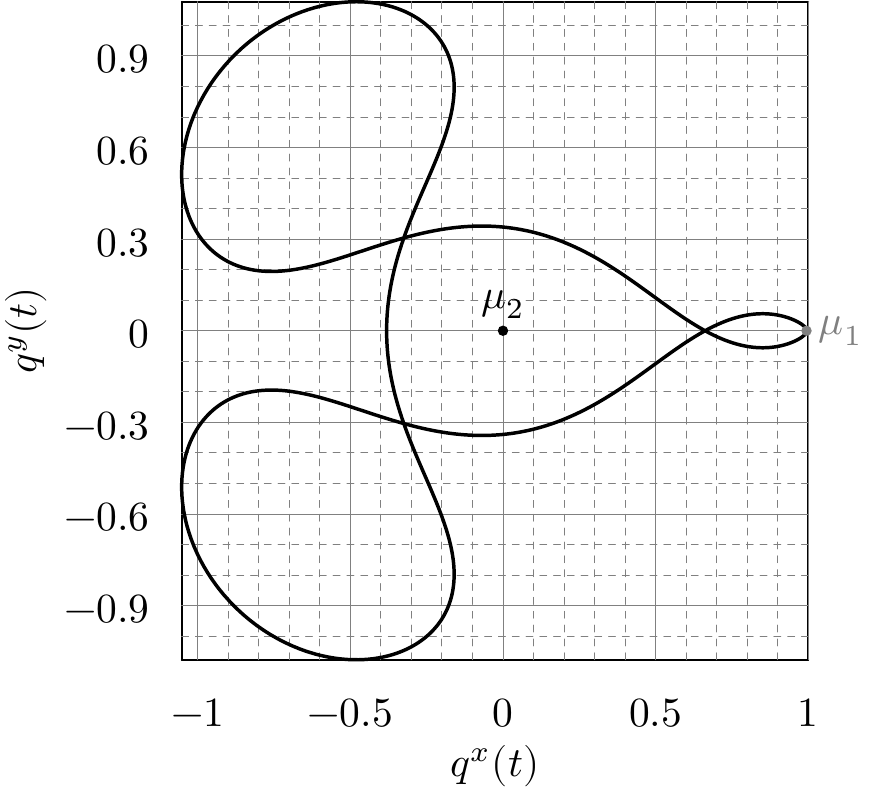}
    \caption{Орбита для второй группы начальных значений~\eqref{eq:arenstorf02:initial}}
    \label{fig:arenstorf02}
  \end{minipage}
  \\
  \begin{minipage}{\linewidth}
    \centering
    \includegraphics[width=0.7\linewidth]{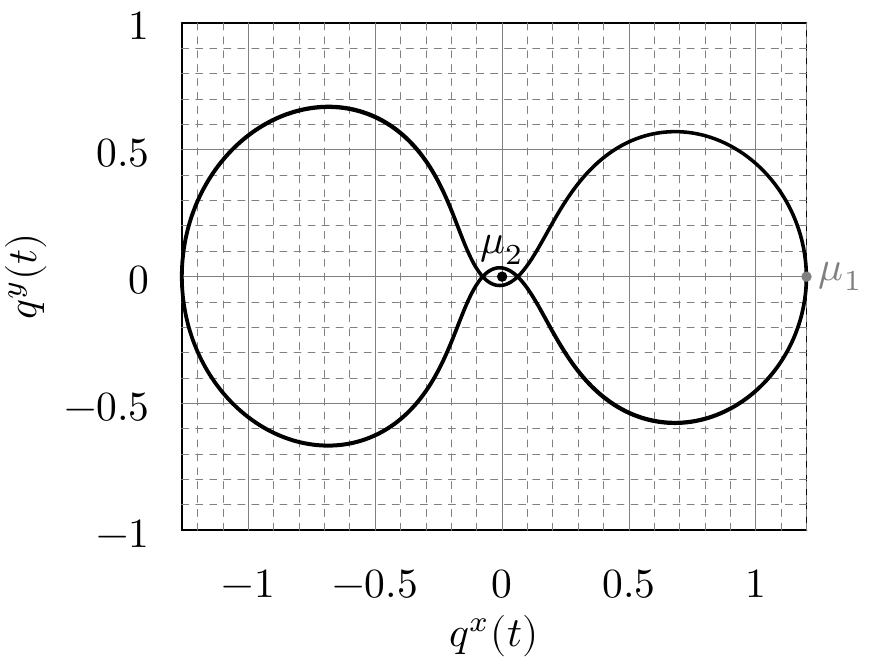}
    \caption{Орбита для третьей группы начальных значений~\eqref{eq:arenstorf03:initial}}
    \label{fig:arenstorf03}
  \end{minipage}
\end{figure}

Численное решение проводится для
$0 \leqslant t \leqslant 17.065216560157962558$ с абсолютной и
относительной погрешностями $A_{tol} = 10^{-17}$ и $R_{tol} = 0$ и
средней массой $\mu_1 = 0.012277471$.
В результате получаем орбиты Аренсторфа
(рис.~\ref{fig:arenstorf01},
\ref{fig:arenstorf02} и \ref{fig:arenstorf03})
для начальных значений~\eqref{eq:arenstorf01:initial},
\eqref{eq:arenstorf02:initial} и \eqref{eq:arenstorf03:initial} соответственно.

Для оценки погрешности система решается численно для временного
промежутка, равного одному периоду. В конечной точке промежутка малое
тело должно вернуться в начальную точку, т.е.
$(q^x(0),q^y(0)) = (q^x(T),q^y(T))$, поэтому погрешность
$\| \mathbf{q}_{0} - \mathbf{q}_{n} \|$ даст глобальную погрешность
метода. Значения погрешности для первой группы начальных значений
приведено в табл.~\ref{tab:errors}.

\begin{table}[!h]
  \caption{Глобальные погрешности численных методов за один оборот по орбите для первой группы начальных значений}
  \label{tab:errors}
  \begin{center}
    \begin{tabular}{ll}
      \hline
      \textbf{Метод} & \textbf{Погрешность}
\\                        
      \hline
      \verb|DPRK546S| & $1.84566 \cdot 10^{-12}$\\
      \verb|DPRK547S| & $2.93634 \cdot 10^{-12}$\\
      \verb|DPRK658M| & $4.22771 \cdot 10^{-12}$\\
      \verb|Fehlberg45| & $7.42775 \cdot 10^{-12}$\\
      \verb|DOPRI5| & $1.95463 \cdot 10^{-13}$\\
      \verb|DVERK65| & $1.67304 \cdot 10^{-12}$\\
      \verb|Fehlberg78B| & $5.45348 \cdot 10^{-12}$\\
      \verb|DOPRI8| & $1.06343 \cdot 10^{-11}$\\
      \hline
    \end{tabular}
  \end{center}
\end{table}
\section{Заключение}

Таким образом, можно сформулировать  основные результата нашей работы:
\begin{enumerate}
\item В качестве специализированной системы символьных вычисление
  применён язык шаблонизатора Jinja2. Относительно использования
  универсальной системы символьных вычислений данное решение позволило
  сделать программный комплекс более простым, компактным и увеличить
  его переносимость.

\item Создан набор сценариев на языке Python с использованием языка
  шаблонизатора Jinja2, которые генерируют код
  на языке Julia, реализующий численные схемы  методов
  Рунге--Кутты без управления шагом, вложенных методов с управлением
  шагом и методов Розенброка также с управлением шагом (все программы
  находятся в открытом доступе и расположены по адресу
  \url{https://bitbucket.org/mngev/rungekutta_generator}).

\item Сгенерированные функции протестированы с помощью решения
  нескольких типовых задач (непосредственно в тексте статьи подробно
  разобрана задача Аренсторфа, как наиболее требовательная к точности
  численной схемы); модуль для языка Julia также доступен по адресу
  \url{https://bitbucket.org/mngev/rungekutta-autogen}.
\end{enumerate}

\def\acknowledgmentsname{Благодарности}

\begin{acknowledgments}

Публикация подготовлена при поддержке Программы РУДН «5-100»
и при финансовой поддержке РФФИ в рамках научного проекта
\textnumero~19-01-00645.

\end{acknowledgments}

  \ifthenelse{\boolean{xetex} \OR \boolean{luatex}}
  {\bibliographystyle{ugost2008l}}
  {\bibliographystyle{gost2008l}}

\bibliography{bib/symbolic-rk/cite}